\documentclass[11pt]{article}
\usepackage{epsfig}

\def\be{\begin{equation}}
\def\ee{\end{equation}}
\def\bea{\begin{eqnarray}}
\def\eea{\end{eqnarray}}
\def\bes{\begin{eqnarray*}}
\def\ees{\end{eqnarray*}}

\def\nn{\nonumber}
\def\<{\langle}
\def\>{\rangle}
\def\lb{\label}
\def\bs{\setminus}
\def\pt{\partial}

\def\R{{\bf R}}

\def\Z{{\bf Z}}
\def\N{{\bf N}}

\def\Q{{\bf Q}}
\def\T{{\bf T}}

\def\aa{{\alpha}}

\def\ep{{\epsilon}}
\def\lm{{\lambda}}
\def\Lm{{\Lambda}}

\def\sg{{\sigma}}

\def\Sg{{\Sigma}}
\def\vf{{\varphi}}

\def\vep{{\varepsilon}}
\def\vth{{\vartheta}}

\def\H{{\cal H}}

\def\T{{\cal T}}

\def\Nn{{\cal N}}

\def\mul{{\rm mul}}

\def\per{{\rm per}}

\def\crit{{\rm crit}}
\def\span{{\rm span}}

\def\td#1{\tilde{#1}}
\def\wtd#1{\widetilde{#1}}
\def\hb{\vrule height0.18cm width0.14cm $\,$}

\title{Resonance Identities for Closed Characteristics on Compact Star-shaped Hypersurfaces in ${\bf R}^{2n}$}
\author{
Hui Liu$^{1,2}$,\thanks{Partially supported by China Postdoctoral Science Foundation No.2013M540512.
E-mail:huiliu@ustc.edu.cn. } \qquad
Yiming Long$^{2}$,\thanks{Partially supported by NSFC (No.11131004), MCME and LPMC of MOE of China, Nankai
University and BCMIIS of Capital Normal University. E-mail: longym@nankai.edu.cn.} \qquad
Wei Wang$^{3}$\thanks{Partially supported by NSFC (No. 11222105), Foundation for the Authors of National Excellent
Doctoral Dissertations of P. R. China No. 201017. E-mail: alexanderweiwang@gmail.com.  } \\ \\
$^{1}$ School of Mathematical Sciences, University of Science and Technology of China, \\ Hefei, Anhui 230026\\
$^{2}$ Chern Institute of Mathematics and LPMC, Nankai University, Tianjin 300071\\
$^{3}$ Key Laboratory of Pure and Applied Mathematics, \\
School of Mathematical Science, Peking University, Beijing 100871\\
People's Republic of China\\}
\date{}

\begin{document}

\maketitle
\begin{abstract}
{\it Resonance relations among periodic orbits on given energy hypersurfaces are very
important for getting deeper understanding of the dynamics of the corresponding Hamiltonian systems.
In this paper, we establish two new resonance identities for closed characteristics on every
compact star-shaped hypersurface $\Sigma$ in ${\bf R}^{2n}$ when the number of geometrically distinct closed
characteristics on $\Sigma$ is finite, which extend those identities established by C. Viterbo in 1989 for
star-shaped hypersurfaces assuming in addition that all the closed characteristics and their iterates
are non-degenerate, and that by W. Wang, X. Hu and Y. Long in 2007 for strictly convex hypersurfaces in
${\bf R}^{2n}$.}
\end{abstract}


\renewcommand{\theequation}{\thesection.\arabic{equation}}
\renewcommand{\thefigure}{\thesection.\arabic{figure}}

\setcounter{equation}{0}
\section{Introduction and main result}

Let $\Sigma$ be a $C^3$ compact hypersurface in $\R^{2n}$ strictly star-shaped with respect to the origin, i.e.,
the tangent hyperplane at any $x\in\Sigma$ does not intersect the origin. We denote the set of all such
hypersurfaces by $\H_{st}(2n)$, and denote by $\H_{con}(2n)$ the subset of $\H_{st}(2n)$ which consists of all
strictly convex hypersurfaces. We consider closed characteristics $(\tau, y)$ on $\Sigma$, which are solutions
of the following problem
\be
\left\{\matrix{\dot{y}=JN_\Sigma(y), \cr
               y(\tau)=y(0), \cr }\right. \lb{1.1}\ee
where $J=\left(\matrix{0 &-I_n\cr
        I_n  & 0\cr}\right)$, $I_n$ is the identity matrix in $\R^n$, $\tau>0$, $N_\Sigma(y)$ is the outward
normal vector of $\Sigma$ at $y$ normalized by the condition $N_\Sigma(y)\cdot y=1$. Here $a\cdot b$ denotes
the standard inner product of $a, b\in\R^{2n}$. A closed characteristic $(\tau, y)$ is {\it prime}, if $\tau$
is the minimal period of $y$. Two closed characteristics $(\tau, y)$ and $(\sigma, z)$ are {\it geometrically
distinct}, if $y(\R)\not= z(\R)$. We denote by $\T(\Sg)$ the set of all geometrically distinct closed
characteristics on $\Sg$. A closed characteristic $(\tau,y)$ is {\it non-degenerate}, if $1$ is a Floquet
multiplier of $y$ of precisely algebraic multiplicity $2$.

The study on closed characteristics in the global sense started in 1978, when the existence of at least one
closed characteristic was first established on any $\Sg\in\H_{st}(2n)$ by P. Rabinowitz in \cite{Rab1}
and on any $\Sg\in\H_{con}(2n)$ by A. Weinstein in \cite{Wei1} independently, since then the existence of
multiple closed characteristics on $\Sg\in\H_{con}(2n)$ has been deeply studied by many mathematicians, for
example, studies in \cite{EkL1}, \cite{EkH1}, \cite{Szu1}, \cite{LoZ1}, \cite{WHL1}, and \cite{Wan1} for
convex hypersurfaces. We refer readers to the survey paper \cite{Lon5} and the recent \cite{Lon6} of Y. Long
for earlier works and references on this subject.

But for the star-shaped hypersurfaces, one difficulty in the study on the star-shaped hypersurfaces is that
the Maslov-type index and mean index of each closed characteristic may be negative. We are only aware of a few
papers about the multiplicity of closed characteristics. In \cite{Gir1} of 1984 and \cite{BLMR} of 1985,
$\;^{\#}\T(\Sg)\ge n$ for $\Sg\in\H_{st}(2n)$ was proved under some pinching conditions. In \cite{Vit1} of 1989,
C. Viterbo proved a generic existence result for infinitely many closed characteristics on star-shaped
hypersurfaces. In \cite{HuL1} of 2002, X. Hu and Y. Long proved that $\;^{\#}\T(\Sg)\ge 2$ for
$\Sg\in \H_{st}(2n)$ on which all the closed characteristics and their iterates are non-degenerate. Recently
$\;^{\#}\T(\Sg)\ge 2$ was proved for every $\Sg\in \H_{st}(4)$ by D. Cristofaro-Gardiner and M. Hutchings
in \cite{CGH1}, and it's different proofs can also be found in \cite{GHHM}, \cite{LLo1} and \cite{GiG1}.

In \cite{Eke1} of 1984, I. Ekeland first discovered some resonance relations of closed characteristics for
$\Sg\in\H_{con}(2n)$, but which are not explicitly given. In \cite{Vit1}, C. Viterbo established two such
identities explicitly for closed characteristics on $\Sg\in\H_{st}(2n)$ under the assumption that all the
closed characteristics on $\Sg$ are non-degenerate. Such identities are important ingredients in the study
in \cite{Vit1} and \cite{HuL1}. In \cite{WHL1} of 2007, W. Wang, X. Hu and Y. Long proved the resonance
identity for every $\Sg\in\H_{con}(2n)$ which removed the non-degeneracy condition. This identity plays a
crucial role in the proof of their estimate $\,^{\#}\T(\Sg)\ge 3$ for every $\Sg\in \H_{con}(6)$. 
Note that in \cite{Rad1} of 1989 and \cite{Rad2} of 1992, a similar identity for closed
geodesics on compact Finsler manifolds was established by H.-B. Rademacher.
Motivated by \cite{Vit1} and \cite{WHL1}, we establish the following new identities on closed characteristics for
star-shaped hypersurface $\Sg\in\H_{st}(2n)$ without the non-degeneracy conditions.

{\bf Theorem 1.1.} {\it Suppose that $\Sg\in \H_{st}(2n)$ satisfies $\,^{\#}\T(\Sg)<+\infty$. Denote all
the geometrically distinct prime closed characteristics by $\{(\tau_j,\; y_j)\}_{1\le j\le k}$. Then the
following identities hold
\bea
\sum_{1\le j\le k\atop \hat{i}(y_j)>0}\frac{\hat{\chi}(y_j)}{\hat{i}(y_j)} &=& \frac{1}{2},  \lb{1.2}\\
\sum_{1\le j\le k\atop \hat{i}(y_j)<0}\frac{\hat{\chi}(y_j)}{\hat{i}(y_j)} &=& 0,  \lb{1.3}\eea
where $\hat{i}(y_j)\in\R$ is the mean index of $y_j$ given by Definition 4.7, $\hat{\chi}(y_j)\in\Q$ is
the average Euler characteristic given by Definition 4.8 and Remark 4.9 below. Especially by (\ref{4.32}),
we have
\be \hat\chi(y) = \frac{1}{K(y)}\sum_{1\le m\le K(y)\atop 0\le l\le 2n-2}(-1)^{i(y^{m})+l}k_l(y^{m}),  \lb{1.4}\ee
where $K(y)\in \N$ is the minimal period of critical modules of iterations of $y$ defined in Proposition 4.6,
$i(y^{m})$ is the index defined in Definition 4.7 (cf. Definition 2.9 and (\ref{2.15}) below), and $k_l(y^{m})$s
are the critical type numbers of $y^m$ given by Definition 4.3 and Remark 4.4 below. }

{\bf Remark 1.2.} When all the closed characteristics on $\Sigma\in\H_{st}(2n)$ together with their iterations
are non-degenerate, by Remark 4.9 our identities (\ref{1.2}) and (\ref{1.3}) coincide with the identities (1.3)
and (1.4) of Theorem 1.2 of \cite{Vit1}. Thus our Theorem 1.1 generalizes C. Viterbo's result in \cite{Vit1}
to the degenerate case.

When $\Sigma\in\H_{con}(2n)$, we can choose $K_0=0$ in the proof of Case (b) of Theorem 3.3 below. Then $d(K)=0$
in (\ref{2.15}). By (\ref{3.13}) and (\ref{3.15}), we obtain
$$  C_{S^1, l}(F_{K},S^1\cdot\bar{x})\cong C_{S^1,l}(\wtd{F}_{K}, S^1\cdot\bar{y})\cong C_{S^1,l}(\wtd{F}_{0}, S^1\cdot\bar{y}).  $$
Noticing that $C_{S^1, l}(\wtd{F}_{0}, S^1\cdot\bar{y})$ is exactly isomorphic to
$C_{S^1, l}(\Psi_a, S^1\cdot \dot{\bar{x}})$ which is defined in Definition 3.1 of \cite{WHL1}, then our identity
(\ref{1.2}) coincides with the identity (1.3) of Theorem 1.2 of \cite{WHL1}. Thus our Theorem 1.1 generalizes also
the resonance identity in \cite{WHL1} for convex hypersurfaces to star-shaped hypersurfaces.

We also note that some similar resonance identities for closed Reeb orbits on closed contact manifolds were
established in Theorem 3.6 of \cite{GiG1} under the context of local contact homology after we completed this paper.

The main idea in the proof of Theorem 1.1 and the arrangement of the rest of this paper are as follows.

$\<1\>$ Motivated by the works \cite{Vit1} of C. Viterbo and \cite{WHL1} of W. Wang, X. Hu, and Y. Long, for
every $\Sg\in\H_{st}(2n)$ with $\,^{\#}\T(\Sg)<+\infty$, we shall construct a functional $F_{a, K}$ on the space
$W^{1,2}(\R/\Z, \R^{2n})$ for large $a>0$ and $K$ satisfying the requirement (\ref{2.3})-(\ref{2.4}) below and
establish a Morse theory of this functional $F_{a, K}$ to study closed characteristics on $\Sg$. By a change of
variable, it is equivalent to study a functional $\Psi_{a, K}$ on $L^{2}(\R/\Z, \R^{2n})$.

As usual we use the Clarke-Ekeland dual action principle and a modification of the Viterbo index
theory. Because in general such a dual action functional is not $C^2$, motivated by the studies on
closed geodesics and convex Hamiltonian systems, we follow \cite{Vit1} to introduce a finite dimensional
approximation to the space $L^2(\R/\Z, \R^{2n})$ to get the enough smoothness. This finite dimensional
approximation allows us to apply the idea of the Splitting Lemma of D. Gromoll and W. Meyer \cite{GrM1}
to obtain the periodicity of critical modules for closed characteristics, which overcomes the first difficulty
in addition to the study in \cite{WHL1}. The second difficulty and the most important thing is that, all the
critical modules at a critical orbit $S^1\cdot x$ of $F_{a,K}$ rely on $K$, and we need to show they are isomorphic
to each other and thus are independent of such $K$. This is proved by Theorem 3.3 below.

Because the functional $F_{a,K}$ is not $C^2$ on $W^{1, 2}(\R/\Z, \R^{2n})$,
we can not get Splitting Lemma for $F_{a,K}$ directly. But in our case, the functional
$F_{a,K^\prime}$ is uniformly concave in the direction of $D_\infty(K_0)$ by (\ref{3.7}) below,
where $D_\infty(K_0)$ is as in Definition 3.4, $K_0<K^\prime$ satisfying that $K^\prime-K_0$ is small enough.
Motivated by the method of \cite{DHK1},  we obtain a Splitting Lemma type argument(cf. Lemma 3.5) to complete the proof of Theorem 3.3.

$\<2\>$ To achieve the above mentioned purposes, following Proposition 2.2 of \cite{WHL1} and Proposition 2.7 of
\cite{Vit1}, we construct a special family of Hamiltonian functions which have more flexible properties at the
origin and infinity, and are homogenous in the middle and near the critical orbits.

In Section 2, fixing a hypersurface $\Sg\in\H_{st}(2n)$ with $\,^{\#}\T(\Sg)<+\infty$, we construct a family
of Hamiltonian functions in Lemma 2.4 using auxiliary functions satisfying conditions (i)-(ii) of Lemma 2.2,
together with Proposition 2.5 which yields more precise requirement on the Hamiltonian functions near the origin
and infinity. Using such modified Hamiltonian functions, we construct a functional $F_{a, K}$ on the space
$W^{1,2}(\R/\Z, \R^{2n})$ for every $a>0$ and $K$ satisfying (\ref{2.3})-(\ref{2.4}), whose critical points are
precisely all the closed characteristics on $\Sg$ with periods less than $aT$ and that the origin of
$W^{1, 2}(\R/\Z, \R^{2n})$ is the only constant critical point of $F_{a, K}$. By a usual change of variables,
properties of $F_{a,K}$ can be studied by a functional $\Psi_{a,K}$ on $L^{2}(\R/\Z, \R^{2n})$. Using the finite
dimensional approximation, we get the Palais-Smale condition for $F_{a,K}$ and prove that for every fixed closed
characteristic $(\tau,y)$ on $\Sg$, the Viterbo index and nullity of all the functionals $F_{a, K}$ at its
critical point corresponding to $(\tau,y)$ are independent of $a$ whenever $a>\frac{\tau}{T}$.

$\<3\>$ In Section 3, we prove that for every fixed closed characteristic $(\tau,y)$ on $\Sg$, the critical
modules of all the functionals $F_{a, K}$ at its critical point corresponding to $(\tau,y)$ are independent of
$a$ and $K$. Here the main difficulty part is to deal with the case when $K$ crosses values in $(2\pi/T)\Z$.
Here the main idea is to use the Splitting Lemma type argument to obtain the independence of critical modules
in $K$.

$\<4\>$ In Section 4, we further require the Hamiltonian function to be homogeneous near every critical orbit
so that the critical modules are periodic functions of the dimension. This homogeneity of the Hamiltonian
function is realized by the condition (iii) of Lemma 2.2.

$\<5\>$ In Section 5, we get a degenerate version of Theorem 7.1 of \cite{Vit1} which shows that the
origin has in fact no homological contribution to the lower order terms in the Morse series.

$\<6\>$ In Section 6, we use the homological information obtained in the Sections 2-5, compute all the local
critical modules of the dual action functional $F_{a, K}$ and use such information to set up a Morse theory
for all the closed characteristics on $\Sg\in\H_{st}(2n)$. Together with the global homological information,
we establish the claimed mean index identities (\ref{1.2})-(\ref{1.3}) and prove Theorem 1.1.

In this paper, let $\N$, $\N_0$, $\Z$, $\Q$, $\R$, and $\R^+$ denote the sets of natural integers, non-negative
integers, integers, rational numbers, real numbers, and positive real numbers respectively. Denote by $(a, b)$
and $|a|$ the standard inner product and norm in $\R^{2n}$. Denote by $\<\cdot,\cdot\>$ and $\|\cdot\|$
the standard $L^2$ inner product and $L^2$ norm. For an $S^1$-space $X$, we denote by $X_{S^1}$ the homotopy
quotient of $X$ by $S^1$, i.e., $X_{S^1}=S^\infty\times_{S^1}X$, where $S^\infty$ is the unit sphere in an
infinite dimensional {\it complex} Hilbert space. By $t\to a^+$, we mean $t>a$ and $t\to a$. In this paper we
use only $\Q$ coefficients for all homological modules.

\setcounter{equation}{0}
\section{Variational structure for closed characteristics and finite dimensional reduction}

In this paper, we follow the frame works of \cite{WHL1} and \cite{Vit1} to transform the problem (\ref{1.1})
into a fixed period problem of a Hamiltonian system with some period $T>0$, which is fixed for the rest of
the paper without further restrictions, and then study its variational structure. Here we omit most of the
details and only point out differences from \cite{WHL1} when necessary.

In the rest of this paper, we fix first a $\Sg\in\H_{st}(2n)$ and assume the following condition on $\T(\Sg)$:

\noindent (F) {\bf There exist only finitely many geometrically distinct prime closed characteristics
\\$\quad \{(\tau_j, y_j)\}_{1\le j\le k}$ on $\Sigma$. }

As in \cite{WHL1}, we have the following discrete subset of $\R^+$:

{\bf Definition 2.1} {\it Under the assumption (F), the set of periods of closed characteristics on $\Sigma$
is defined by
$$ \per(\Sigma)=\{m\tau_j\;|\; m\in\N,\; 1\le j\le k\}. $$
and let $\hat{\tau}=\inf\{s\,|\, s\in \per(\Sigma)\}$.}

Motivated by Lemma 2.2, Proposition 2.7 of \cite{Vit1}, and omitting the condition (iv) of Proposition 2.2 of
\cite{WHL1} to get more flexibility, we use the following auxiliary function to further define Hamiltonian
functions.

{\bf Lemma 2.2} {\it For any sufficiently small $\vth\in (0,1)$, there exists a function
$\vf\equiv \vf_{\vth}\in C^\infty(\R, \R^+)$ depending on $\vth$ which has $0$ as its
unique critical point in $[0, +\infty)$ such that the following hold.

(i) $\vf(0)=0=\vf^\prime(0)$, and $\vf^{\prime\prime}(0)=1=\lim_{t\rightarrow 0^+}\frac{\vf^\prime(t)}{t}$;

(ii) $\frac{d}{dt}\left(\frac{\vf^\prime(t)}{t}\right)<0$ for $t>0$, and
$\lim_{t\rightarrow +\infty}\frac{\vf^\prime(t)}{t}<\vth$; that is, $\frac{\vf^\prime(t)}{t}$ is
strictly decreasing for $t> 0$;

(iii) In particular, we can choose $\aa\in (1,2)$ sufficiently close to $2$ and $c\in (0,1)$ such that
$\vf(t)=ct^{\aa}$ whenever $\frac{\vf'(t)}{t}\in [\vth, 1-\vth]$ and $t>0$.}

{\bf Remark 2.3.} As in \cite{WHL1}, Lemma 2.2 (iii) above is used only in our study in the Section 4 to obtain
the periodic property of critical modules at critical points. In the other parts of this paper we use function
$\vf$ which satisfy the properties (i)-(ii) only and defined on $[0,\,+\infty)$. In the proof of Lemma 2.4, given
an $a>\frac{\hat{\tau}}{T}$, we choose first the parameter $\vth\in (0,\frac{\hat{\tau}}{aT})$ depending on $a$.
Then we choose the parameter $\alpha\in (1,\,2)$ depending on $a$ so that the proof of Lemma 2.2 goes through,
and choose $\vf$ to be homogeneous of degree $\alpha$ and modify it near $0$ and $+\infty$ such that (i)-(ii) of
Lemma 2.2 hold. We denote such choices of $\vth$, $\aa$ and $\vf$ by $\vth_a$, $\aa_a$ and $\vf_a$ respectively to
indicate their dependence on $a$. In such a way, we can obtain a connected family of $\vf_a$ satisfying (i)-(ii)
of Lemma 2.2 such that $\vf_a$ and its first and second derivatives with respect to $t$ depend continuously on $a$.

Let $j: \R^{2n}\rightarrow\R$ be the gauge function of $\Sigma$, i.e., $j(\lambda x)=\lambda$ for $x\in\Sigma$
and $\lambda\ge0$, then $j\in C^3(\R^{2n}\bs\{0\}, \R)\cap C^0(\R^{2n}, \R)$ and $\Sigma=j^{-1}(1)$. Then
the following lemma was proved in Proposition 2.4 (iii) of \cite{WHL1} (cf. also Lemmas 2.1 and 2.2 of \cite{Vit1}).

{\bf Lemma 2.4.} {\it Let $a>\frac{\hat{\tau}}{T}$, $\vth_a\in (0, \frac{\hat{\tau}}{aT})$ and $\vf_a$ be a
$C^\infty$ function associated to $\vth_a$ satisfying (i)-(ii) of Lemma 2.2 and continuously depending on
the parameter $a$ as mentioned in Remark 2.3. Define the Hamiltonian function $\wtd{H}_a(x)=a\vf_a(j(x))$
and consider the fixed period system
\be \left\{\matrix{\dot{x}(t) &=& J\wtd{H}_a^\prime(x(t)), \cr
                   x(0) &=& x(T).\qquad    \cr }\right.   \lb{2.1}\ee
Then solutions of (\ref{2.1}) are $x\equiv 0$ and $x=\rho y(\tau t/T)$ with
$\frac{\vf_a^\prime(\rho)}{\rho}=\frac{\tau}{aT}$, where $(\tau, y)$ is a solution of (\ref{1.1}). In particular,
non-zero solutions of (\ref{2.1}) are in one to one correspondence with solutions of (\ref{1.1}) with period
$\tau<aT$.}

For technical reasons we want to further modify the Hamiltonian, more precisely, we follow Page 624 of \cite{Vit1},
and let $\ep_a$ satisfy $\ep_a T<2\pi$ and $\vth_a$ be small enough, we can construct a function $H_a$, coinciding
with $\wtd{H}_a$ on $U_A=\{x\mid \wtd{H}_a(x)\leq A\}$ for some large $A$, and with $\frac{1}{2}\ep_a|x|^2$ outside
some large ball, such that $\nabla H_a(x)$ does not vanish and $H_a^{\prime\prime}(x)<\ep_a$ outside $U_A$. As in
Proposition 2.7 of \cite{Vit1}, we have the following result.

{\bf Proposition 2.5.} {\it For $a>\frac{\hat{\tau}}{T}$ and small $\ep_a$, we choose small enough $\vth_a$ such
that Lemma 2.4 holds. Then there exists a function $H_a$ on ${\bf R}^{2n}$ such that $H_a$ is $C^1$ on ${\bf R}^{2n}$,
and $C^3$ on ${\bf R}^{2n}\bs\{0\}$, $H_a=\wtd{H}_a$ in $U_A$, and $H_a(x)=\frac{1}{2}\ep_a|x|^2$ for $|x|$ large,
and the solutions of the fixed period system
\be   \left\{\matrix{\dot{x}(t) &=& JH_a^\prime(x(t)), \cr
                   x(0) &=& x(T),\qquad         \cr }\right.  \lb{2.2}\ee
are the same with those of (\ref{2.1}).}

{\bf Remark 2.6.} Note that here the first derivative of $H_a(x)$ with respect to $x\in\R^{2n}$ and
the second derivative of $H_a(x)$ with respect to $x\in\R^{2n}\bs\{0\}$ are continuous in the parameter
$a$. Note that under these choices, the first and second derivatives of $\ep_a$ with respect to $a$ are also
continuous. Here, that $H_a$'s form a connected family in $a$ is crucial in our proofs below for Lemma 2.10,
and Propositions 2.11 and 3.2.

As in \cite{BLMR} (cf. Section 3 of \cite{Vit1}), for any $a>\frac{\hat{\tau}}{T}$, we can choose some
large constant $K=K(a)$ such that
\be H_{a,K}(x) = H_a(x)+\frac{1}{2}K|x|^2   \lb{2.3}\ee
is a strictly convex function, that is,
\be (\nabla H_{a, K}(x)-\nabla H_{a, K}(y), x-y) \geq \frac{\ep}{2}|x-y|^2,  \lb{2.4}\ee
for all $x, y\in {\bf R}^{2n}$, and some positive $\ep$. Let $H_{a,K}^*$ be the Fenchel dual of $H_{a,K}$
defined by
\be  H_{a,K}^\ast (y) = \sup\{x\cdot y-H_{a,K}(x)\;|\; x\in \R^{2n}\}.   \lb{2.5}\ee
The dual action functional on $X=W^{1, 2}({\bf R}/{T {\bf Z}}, {\bf R}^{2n})$ is defined by
\be F_{a,K}(x) = \int_0^T{\left[\frac{1}{2}(J\dot{x}-K x,x)+H_{a,K}^*(-J\dot{x}+K x)\right]dt}.  \lb{2.6}\ee
Then $F_{a,K}\in C^{1,1}(X,{\bf R})$ holds as proved in (3.16) of \cite{Vit1}, but $F_{a,K}$ is not $C^2$.

{\bf Lemma 2.7.} (cf. Proposition 3.4 of \cite{Vit1}) {\it Assume $KT\not\in 2\pi{\bf Z}$, then $x$ is a
critical point of $F_{a, K}$ if and only if it is a solution of (\ref{2.2}).}

From Lemma 2.7, we know that the critical points of $F_{a,K}$ are independent of $K$.

{\bf Proposition 2.8.} {\it For every critical point $x_a\neq 0$ of $F_{a,K}$, the critical value
$F_{a,K}(x_a)<0$ holds and is independent of $K$.}

{\bf Proof.} Since $\nabla H_{a,K}(x_a)=-J\dot{x}_a+Kx_a$, then we have
$$  H_{a,K}^*(-J\dot{x}_a+Kx_a)=(-J\dot{x}_a+Kx_a, x_a)-H_{a,K}(x_a).  $$
Thus we obtain
\bea F_{a,K}(x_a)
&=& \int_0^T{\left[\frac{1}{2}(J\dot{x}_a-Kx_a,x_a)+H_{a,K}^*(-J\dot{x}_a+Kx_a)\right]dt}  \nn\\
&=& \int_0^T{\left[-\frac{1}{2}(J\dot{x}_a-Kx_a,x_a)-H_{a,K}(x_a)\right]dt}  \nn\\
&=& \int_0^T{\left[-\frac{1}{2}(J\dot{x}_a,x_a)-H(x_a)\right]dt}  \nn\\
&=& \int_0^T{\left[\frac{1}{2}(H_a^\prime(x_a),x_a)-H(x_a)\right]dt}.  \lb{2.7}\eea
By Lemma 2.4 and Proposition 2.5, we have $x_a=\rho_a y(\tau t/T)$ with
$\frac{\vf_a^\prime(\rho_a)}{\rho_a}=\frac{\tau}{aT}$. Hence, we have
\be  F_{a,K}(x_a)=\frac{1}{2}a\vf_a^\prime(\rho_a)\rho_aT-a\vf_a(\rho_a)T.   \lb{2.8}\ee
Here we used the facts that $j^\prime(y)=N_\Sigma(y)$ and $j^\prime(y)\cdot y=1$.

Let $f(t)=\frac{1}{2}a\vf_a^\prime(t)t-a\vf_a(t)$ for $t\ge 0$. Then we have $f(0)=0$ and
$f'(t)=\frac{a}{2}(\vf_a^{\prime\prime}(t)t-\vf_a^\prime(t))<0$ since $\frac{d}{dt}(\frac{\vf_a^\prime(t)}{t})<0$
by (ii) of Lemma 2.2. Together with (\ref{2.8}), it yields the proposition. \hfill\hb

As well known, when $KT\notin 2\pi{\bf Z}$, the map $x\mapsto -J\dot{x}+Kx$ is a Hilbert space isomorphism between
$X=W^{1, 2}({\bf R}/{T {\bf Z}}; {\bf R}^{2n})$ and $E=L^{2}({\bf R}/(T {\bf Z}),{\bf R}^{2n})$. We denote its inverse
by $M_K$ and the functional
\be \Psi_{a,K}(u)=\int_0^T{\left[-\frac{1}{2}(M_{K}u, u)+H_{a,K}^*(u)\right]dt}, \qquad \forall\,u\in E. \lb{2.9}\ee
Then $x\in X$ is a critical point of $F_{a,K}$ if and only if $u=-J\dot{x}+Kx$ is a critical point of $\Psi_{a, K}$.
We have a natural $S^1$-action on $X$ or $E$ defined by
\be  \theta\cdot u(t)=u(\theta+t),\quad\forall\, \theta\in S^1, \, t\in\R.  \lb{2.10}\ee
Clearly both of $F_{a, K}$ and $\Psi_{a, K}$ are $S^1$-invariant. For any $\kappa\in\R$, we denote by
\bea
\Lambda_{a, K}^\kappa &=& \{u\in L^{2}({\bf R}/{T {\bf Z}}; {\bf R}^{2n})\;|\;\Psi_{a,K}(u)\le\kappa\}  \lb{2.11}\\
X_{a, K}^\kappa &=& \{x\in W^{1, 2}({\bf R}/(T {\bf Z}),{\bf R}^{2n})\;|\;F_{a, K}(x)\le\kappa\}.  \lb{2.12}\eea
Clearly, both level sets are $S^1$-invariant.

{\bf Definition 2.9.} (cf. p.628 of \cite{Vit1}) {\it Suppose $u$ is a nonzero critical point of $\Psi_{a, K}$.
Then the formal Hessian of $\Psi_{a, K}$ at $u$ is defined by
\be Q_{a,K}(v)=\int_0^T(-M_K v\cdot v+H_{a,K}^{*\prime\prime}(u)v\cdot v)dt,  \lb{2.13}\ee
which defines an orthogonal splitting $E=E_-\oplus E_0\oplus E_+$ of $E$ into negative, zero and positive subspaces.
The index and nullity of $u$ are defined by $i_K(u)=\dim E_-$ and $\nu_K(u)=\dim E_0$ respectively. }

Similarly, we define the index and nullity of $x=M_Ku$ for $F_{a, K}$, we denote them by $i_K(x)$ and
$\nu_K(x)$. Then we have
\be  i_K(u)=i_K(x),\quad \nu_K(u)=\nu_K(x),  \lb{2.14}\ee
which follow from the definitions (\ref{2.6}) and (\ref{2.9}). The following important formula was proved in
Lemma 6.4 of \cite{Vit1}:
\be  i_K(x) = 2n([KT/{2\pi}]+1)+i^v(x) \equiv d(K)+i^v(x),   \lb{2.15}\ee
where the index $i^v(x)$ does not depend on K, but only on $H_a$.

By the proof of Proposition 2 of \cite{Vit2}, we have that $v\in E$ belongs to the null space of $Q_{a, K}$
if and only if $z=M_K v$ is a solution of the linearized system
\be  \dot{z}(t) = JH_a''(x(t))z(t).  \lb{2.16}\ee
Thus the nullity in (\ref{2.14}) is independent of $K$, which we denote by $\nu^v(x)\equiv \nu_K(u)= \nu_K(x)$.

In this paper, we say that $\Psi_{a,K}$ with $a\in [a_1, a_2]$ form a continuous family of functionals
in the sense of Remark 2.6, when $0<a_1<a_2<+\infty$.

Motivated by Proposition 3.9 and Lemma 5.1 of \cite{Vit1} as well as Lemma 3.4 of \cite{WHL1}, we
have the following

{\bf Lemma 2.10.} {\it For any $0<a_1<a_2<+\infty$, let $K$ be fixed so that $\Psi_{a, K}$ with $a\in [a_1, a_2]$
is a continuous family of functionals defined by (\ref{2.9}) satisfying (\ref{2.4}) with the same $\ep>0$. Then
there exist a finite dimensional $S^1$-invariant subspace $G$ of $L^{2}({\bf R}/{T {\bf Z}}; {\bf R}^{2n})$ and
a family of $S^1$-equivariant maps $h_{a}: G\rightarrow G^\perp$ such that the following hold.

(i) For $g\in G$, each function $h\mapsto\Psi_{a,K}(g+h)$ has $h_a(g)$ as the unique minimum in $G^\perp$.

Let $\psi_{a,K}(g)=\Psi_{a,K}(g+h_a(g))$. Then we have

(ii) Each $\psi_{a, K}$ is $C^1$ and $S^1$-invariant on $G$. Here $g_a$ is a critical point of $\psi_{a, K}$
if and only if $g_a+h_{a}(g_a)$ is a critical point of $\Psi_{a, K}$.

(iii) If $g_a\in G$ and $H_a$ is $C^k$ with $k\ge 2$ in a neighborhood of the trajectory of $g_a+h_{a}(g_a)$,
then $\psi_{a,K}$ is $C^{k-1}$ in a neighborhood of $g_a$. In particular, if $g_a$ is a nonzero critical point
of $\psi_{a,K}$, then $\psi_{a,K}$ is $C^2$ in a neighborhood of the critical orbit $S^1\cdot g_a$. The
index and nullity of $\Psi_{a,K}$ at $g_a+h_{a}(g_a)$ defined in Definition 2.9 coincide with the Morse index
and nullity of $\psi_{a,K}$ at $g_a$.

(iv) For any $\kappa\in\R$, we denote by¡¤¡¤
\be  \wtd{\Lambda}_{a,K}^\kappa=\{g\in G \;|\; \psi_{a,K}(g)\le\kappa\}.   \lb{2.17}\ee
Then the natural embedding $\wtd{\Lambda}_{a, K}^\kappa \hookrightarrow {\Lambda}_{a, K}^\kappa $ given by
$g\mapsto g+h_{a}(g)$ is an $S^1$-equivariant homotopy equivalence.

(v) The functionals $a\mapsto\psi_{a,K}$ is continuous in $a$ in the $C^1$ topology. Moreover
$a\mapsto\psi^{\prime\prime}_{a,K}$ is continuous in a neighborhood of the critical orbit $S^1\cdot g_a$.}

{\bf Proof.} Firstly, we consider the eigenvalues of $-M_K$. Let $x(t)=e^{-JLt}x_0$ for some $L\in\frac{2\pi}{T}\Z$
and $x_0\in\R^{2n}$, then $-J\dot{x}+Lx=(L+K)x$. Thus $\{-\frac{1}{L+K}\mid L\in\frac{2\pi}{T}\Z\}$ is the set
of all the eigenvalues of $-M_K$.

By the convexity of $H^*_{a,K}$, we have
\be (H_{a,K}^{*\prime}(u)-H_{a,K}^{*\prime}(v), u-v)\ge
            \omega|u-v|^2,\quad \forall\;a\in[a_1, a_2],\; u, v\in\R^{2n}, \lb{2.18}\ee
for some $\omega>0$. Hence we can use the proof of Proposition 3.9 of \cite{Vit1} to obtain the subspace $G$ and
the map $h_a$. In fact, Let $G$ be the subspace of $L^{2}({\bf R}/(T {\bf Z}); {\bf R}^{2n})$ generated by the
eigenvectors of $-M_K$ whose eigenvalues are less than $-\frac{\omega}{2}$, i.e.,
$$ G = \span\{e^{-JLt}x_0\mid -\frac{1}{L+K}<-\frac{\omega}{2}, L\in\frac{2\pi}{T}\Z, x_0\in\R^{2n}\}, $$
and $h_a(g)$ is defined by the equation
\be \frac{\partial}{\partial h}\Psi_{a, K}(g+h_a(g))=0.  \lb{2.19}\ee
Then (i)-(iii) follows from Proposition 3.9 of \cite{Vit1}, and (iv) follows from Lemma 5.1 of \cite{Vit1}.

The rest part of this proof is devoted to (v).

{\bf Claim. } {\it For each $a\in [a_1, a_2]$ and $\ep>0$ small, we have
\bea
|H^*_{a+\ep, K}(y)-H^*_{a, K}(y)| &=& O(\ep)+O(\ep)|y|^2,\quad \forall y\in\R^{2n}, \lb{2.20}\\
|H^{*\prime}_{a+\ep, K}(y)-H^{*\prime}_{a, K}(y)| &=& O(\ep)+O(\ep)|y|, \quad \forall y\in\R^{2n}, \lb{2.21}\eea
where we denote by $B=O(\ep)$ if $|B|\le C|\ep|$ for some constant $C>0$.}

In fact, we fix an $a\in [a_1, a_2]$ and let $b\in [a-\ep, a+\ep]\cap [a_1, a_2]$. For any $y\in\R^{2n}$,
let $H^{*\prime}_{b, K}(y)=z_b$, then $H^{\prime}_{b, K}(z_b)=y$. By the convexity of $H_{b, K}$, we have
\be |u_1-u_2|\le \alpha |H^{\prime}_{b,K}(u_1)-H^{\prime}_{b, K}(u_2)|, \qquad \forall\,u_1, u_2\in {\bf R}^{2n},
             b\in [a-\ep, a+\ep], \lb{2.22}\ee
for some constant $\alpha>0$ which is independent of $b$. Thus, we obtain
\bea |H^{*\prime}_{a+\ep, K}(y)-H^{*\prime}_{a,K}(y)|
&=& |z_{a+\ep}-z_a|  \nn\\
&\le& \alpha|H^{\prime}_{a+\ep,K}(z_{a+\ep})-H^{\prime}_{a+\ep,K}(z_a)| \nn\\
&=& \alpha|H^{\prime}_{a, K}(z_{a})-H^{\prime}_{a+\ep,K}(z_a)| \nn\\
&=& \alpha(O(\ep)+O(\ep)|z_a|)  \nn\\
&=& \alpha(O(\ep)+O(\ep)\alpha|y|).  \lb{2.23}\eea
Here we have used the fact that $H_a=\frac{1}{2}\ep_a|x|^2$ for $|x|$ large and the derivative of $\ep_a$
with respect to $a$ is continuous by Remark 2.6. Hence, (\ref{2.21}) holds.

For (\ref{2.20}), we have
\be H^{*}_{b, K}(y)=z_b\cdot y-H_{b, K}(z_b).  \lb{2.24}\ee
Then
\bea |H^{*}_{a+\ep, K}(y)-H^{*}_{a,K}(y)|
&=& |(y,z_{a+\ep}-z_a)+ H_{a+\ep,K}(z_{a+\ep})-H_{a,K}(z_a)|   \nn\\
&=& O(\ep)+O(\ep)|y|^2.  \lb{2.25}\eea
Here we used (\ref{2.23}) and the fact that $H_a=\frac{1}{2}\ep_a|x|^2$ for $|x|$ large and
the derivative of $\ep_a$ with respect to $a$ is continuous by Remark 2.6. The claim is proved.

Now we have the following estimates:
\bea
|\Psi_{a+\ep,K}(u)-\Psi_{a,K}(u)| &\le& \int_0^T|H^*_{a+\ep, K}(u)-H^*_{a, K}(u)|dt
     = O(\ep)+O(\ep)\|u\|^2,  \lb{2.26}\\
\|\Psi_{a+\ep,K}^\prime(u)-\Psi_{a,K}^\prime(u)\|^2 &=& \|H^{*\prime}_{a+\ep,K}(u)-H^{*\prime}_{a,K}(u)\|^2
     = O(\ep)+O(\ep)\|u\|^2.  \lb{2.27}\eea
As in \cite{Vit1}, (\ref{2.19}) and the definition of $G$ yield
\be \langle\Psi^\prime_{a,K}(u)-\Psi^\prime_{a, K}(v),\; u-v\rangle
     \ge \frac{\omega}{2}\|u-v\|^2, \quad\forall u-v\in G^\perp,\; a\in [a_1,a_2]. \lb{2.28}\ee
Hence we have
\bea  \frac{\omega}{2}\|h_{a+\ep}(g)-h_a(g)\|^2
&\le& \langle\Psi^\prime_{a+\ep,K}(g+h_{a+\ep}(g))-\Psi^\prime_{a+\ep,K}(g+h_a(g)),\; h_{a+\ep}(g)-h_a(g)\rangle   \nn\\
&=& \langle\Psi^\prime_{a,K}(g+h_a(g))- \Psi^\prime_{a+\ep,K}(g+h_a(g)) ,\; h_{a+\ep}(g)-h_a(g)\rangle  \nn\\
&=& (O(\ep)+O(\ep)\|g+h_a(g)\|^2)^{1/2}\|h_{a+\ep}(g)-h_a(g)\|,   \nn\eea
where the first equality follows from (\ref{2.19}) and the last equality follows from (\ref{2.27}). Hence
the map $a\mapsto h_a(g)$ is continuous.

Because $\psi_{a, K}(g)=\Psi_{a, K}(g+h_a(g))$ by definition,
$\psi_{a, K}^\prime(g)=\frac{\partial}{\partial g}\Psi_{a, K}(g+h_a(g))$ by (\ref{2.19}), hence the first
statement of (v) follows from (\ref{2.26}) and (\ref{2.27}). The last statement of (v) follows from p.629
of \cite{Vit1} and the implicit functional theorem with parameters. \hfill\hb

{\bf Proposition 2.11.} {\it For all $b\ge a>\frac{\tau}{T}$, let $F_{b,K}$ be the functional defined by
(\ref{2.6}), and $x_b$ be the critical point of $F_{b,K}$ so that $x_b$ corresponds to a fixed closed
characteristic $(\tau,y)$ on $\Sigma$ for all $b\ge a$. Then the index $i^v(x_b)$ and nullity $\nu^v(x_b)$
are constants for all $b\ge a$. In particular, when $H_b$ is $\aa$-homogenous for some $\aa\in (1,2)$ near
the image set of $x_b$, the index and nullity coincide with those defined for the Hamiltonian
$H(x)=j(x)^\alpha$ for all $x\in\R^{2n}$. Especially $1\le \nu^v(x_b)\le 2n-1$ always holds.}

{\bf Proof.} Denote by $R(t)$ the fundamental solution of (\ref{2.16}) satisfying $R(0)=I_{2n}$. Then by
Lemma 1.6.11 of \cite{Eke2}, whose proof does not need the convexity of $\Sigma$, we have
\be  R(t)T_{y(0)}\Sigma\subset T_{y(\tau t)}\Sigma.   \lb{2.29}\ee
Then the completely same argument of Proposition 3.5 of \cite{WHL1} proves that $\nu^v(x_a)$ is constant
for all $H_a$ satisfying Proposition 2.5 with $a>\frac{\tau}{T}$ and $1\le \nu^v(x_a)\le 2n-1$.

For any $b>a>\frac{\tau}{T}$, by (iii) of Lemma 2.2, we can construct a continuous family of $\Psi_{c, K}$ with
$c\in[a, b]$ such that $H_b$ is homogenous of degree $\aa=\aa_b$ near the image set of $x_b$. Now we can use
Lemma 2.10 (v) to obtain a continuous family of $\psi_{c,K}$ such that $\psi_{c,K}^{\prime\prime}(g_c)$ depends
continuously on $c\in [a,b]$, where $g_c$ is the critical point of $\psi_{c,K}$ corresponding to $M_K^{-1}{x}_c$.
Because $\dim\ker\psi_{c,K}^{\prime\prime}(g_c)=\nu_K(M_K^{-1}{x}_c)=\nu^v({x_c})={\rm constant}$, the index of
$\psi_{c,K}^{\prime\prime}(g_c)=i_K(M_K^{-1}{x}_c)=i^v({x_c})+d(K)$ must be constant too. Thus $i^v(x_b)$ is
constant for all $b\ge a$. Note that here we used (\ref{2.14}), (\ref{2.15}), the definition of $i^v(x_b)$
below (\ref{2.16}), and Lemma 2.10 (iii). Since the index $i^v(x_b)$ and nullity $\nu^v(x_b)$ only depend on the
value of $H_b$ near the image set of $x_b$ (cf. Proposition 2 of \cite{Vit2}), then the index and nullity
coincide with those defined for the Hamiltonian $H(x)=j(x)^\alpha$, $\forall~x\in\R^{2n}$. The proof is complete.
\hfill\hb

{\bf Proposition 2.12.} {\it $\Psi_{a, K}$ satisfies the Palais-Smale condition on $E$, and $F_{a, K}$ satisfies
the Palais-Smale condition on $X$, when $KT\notin 2\pi{\bf Z}$.}

{\bf Proof.} We first prove that $\Psi_{a, K}$ satisfies the Palais-Smale condition on
$E=L^{2}({\bf R}/(T {\bf Z}); {\bf R}^{2n})$. Below we use short hand notations $\Psi$, $\psi$, $h$, $H_K^*$,
$F$, and $M$ for $\Psi_{a, K}$, $\psi_{a, K}$, $h_a$, $H_{a, K}^*$, $F_{a, K}$, and $M_K$ respectively.

Assume that $x_j=g_j+h_j\in E$ is a sequence such that $\Psi^\prime(x_j)\rightarrow 0$, as $j\to\infty$, where
$g_j\in G$, $h_j\in G^\bot$. Then
\be \<\Psi^\prime(x_j),h_j-h(g_j) \> = o(1)\|h_j-h(g_j)\|,   \lb{2.30}\ee
where we denote by $B_j=o(1)$ if $B_j\to 0$ as $j\to \infty$. On the other hand, since $h_j-h(g_j)\in G^\bot$
and $\frac{\partial\Psi}{\partial h}(g_j+h(g_j))=0$, we obtain
\be \<\Psi^\prime(x_j), h_j-h(g_j)\> = \<\Psi^\prime(x_j)-\Psi^\prime(g_j+h(g_j)), h_j-h(g_j)\>. \lb{2.31}\ee
From (\ref{2.28}) which implies the convexity of $\Psi$ in the direction of $G^\bot$, we have
\be \<\Psi^\prime(x_j)-\Psi^\prime(g_j+h(g_j)), h_j-h(g_j)\> \ge \frac{\omega}{2}\| h_j-h(g_j)\|^2, \lb{2.32}\ee
for some $\omega>0$. Combining (\ref{2.30})-(\ref{2.32}), we obtain
\be   \|h_j-h(g_j)\| = o(1).  \lb{2.33}\ee
Similar to (3.16) of \cite{Vit1}, because $\nabla H_K^*$ is Lipschitz, we have
\be \|\Psi^\prime(x_j)-\Psi^\prime(g_j+h(g_j))\|\leq C\|h_j-h(g_j)\|.  \lb{2.34}\ee
Now $\Psi^\prime(x_j)\rightarrow 0$, from (\ref{2.33}) and (\ref{2.34}), we have
$\psi^\prime(g_j)=\Psi^\prime(g_j+h(g_j))=o(1)$. Together with Proposition 4.1 of \cite{Vit1}, whose proof
goes through in our setting without any modifications, it shows that $g_j$ has a converging subsequence.
Hence, $x_j=g_j+h_j$ must have a converging subsequence by (\ref{2.33}), i.e., $\Psi$ satisfies the
Palais-Smale condition on $E$.

For the Palais-Smale condition of $F$ on $X=W^{1, 2}({\bf R}/(T {\bf Z}); {\bf R}^{2n})$, we have first
\be  F(Mu)=\Psi(u), \qquad \forall\;u\in E,  \lb{2.35}\ee
where $M$ is a Hilbert space isomorphism between $E$ and $X$. Then
\be \< F^\prime(Mu), Mv\>_X = \<\Psi^\prime(u), v\>_E,  \qquad \forall\;u, v \in E. \lb{2.36}\ee
Now we use the standard $L^2$-norm for $E$, and the norm $\|M^{-1}x\|_E$ for $x\in X$ which is equivalent
to the standard one. Then $\<F^\prime(Mu), Mv\>_X =\<M^{-1}F^\prime(Mu), M^{-1}Mv\>_E$ holds. Together with
(\ref{2.35}) it yields the following identity on $E$:
\be  M^{-1}F^\prime(Mu) = \Psi^\prime(u).  \lb{2.37}\ee
Because $M$ is a Hilbert space isomorphism between $E$ and $X$, by (\ref{2.35}) and (\ref{2.37}) the
Palais-Smale condition of $F$ on $X$ follows from that of $\Psi$ on $E$. \hfill\hb

\setcounter{equation}{0}
\section{Parameter independence of critical modules for closed characteristics}

For a critical point $u$ of $\Psi_{a, K}$ and the corresponding $x=M_K u$ of $F_{a, K}$, let
\bea
\Lm_{a,K}(u) &=& \Lm_{a,K}^{\Psi_{a, K}(u)}
   = \{w\in L^{2}(\R/(T\Z), \R^{2n}) \;|\; \Psi_{a, K}(w)\le\Psi_{a,K}(u)\},  \lb{3.1}\\
X_{a,K}(x) &=& X_{a,K}^{F_{a,K}(x)} = \{y\in W^{1, 2}(\R/(T\Z), \R^{2n}) \;|\; F_{a,K}(y)\le F_{a,K}(x)\}. \lb{3.2}\eea
Clearly, both sets are $S^1$-invariant. Denote by $\crit(\Psi_{a, K})$ the set of critical points of $\Psi_{a, K}$.
Because $\Psi_{a,K}$ is $S^1$-invariant, $S^1\cdot u$ becomes a critical orbit if $u\in \crit(\Psi_{a, K})$.
Note that by the condition (F), Lemma 2.4, Proposition 2.5 and Lemma 2.7, the number of critical orbits of $\Psi_{a, K}$
is finite. Hence as usual we can make the following definition.

{\bf Definition 3.1.} {\it Suppose $u$ is a nonzero critical point of $\Psi_{a, K}$, and $\Nn$ is an $S^1$-invariant open
neighborhood of $S^1\cdot u$ such that $\crit(\Psi_{a,K})\cap (\Lm_{a,K}(u)\cap \Nn) = S^1\cdot u$. Then the $S^1$-critical
modules of $S^1\cdot u$ is defined by
\bea C_{S^1,\; q}(\Psi_{a, K}, \;S^1\cdot u)
&=& H_{S^1,\; q}(\Lambda_{a, K}(u)\cap\Nn,\; (\Lambda_{a, K}(u)\setminus S^1\cdot u)\cap \Nn)\nn\\
&\equiv& H_{q}((\Lambda_{a, K}(u)\cap\Nn)_{S^1},\; ((\Lambda_{a,K}(u)\setminus S^1\cdot u)\cap\Nn)_{S^1}), \lb{3.3}\eea
where $H_{S^1,\;\ast}$ is the $S^1$-equivariant homology with rational coefficients in the sense of A. Borel (cf. Chapter
IV of \cite{Bor1}). Similarly, we define the $S^1$-critical modules $C_{S^1,\; q}(F_{a, K}, \;S^1\cdot x)$ of $S^1\cdot x$
for $F_{a, K}$.}

As well-known, this definition is independent of the choice of $\Nn$ by the excision property of the singular homology
theory (cf. Definition 1.7.5 of \cite{Cha1}). Recall that $X_{S^1}$ is defined at the end of Section 1.

We have the following for critical modules.

{\bf Proposition 3.2.} {\it Let $(\tau,y)$ be a closed characteristic on $\Sigma$. For any $\frac{\tau}{T}<a_1<a_2<+\infty$,
let $K$ be a fixed sufficiently large real number so that (\ref{2.4}) holds for all $a\in [a_1, a_2]$. Then the critical
module $C_{S^1,\; q}(F_{a, K}, \;S^1\cdot x)$ is independent of the choice of $H_a$ defined in Proposition 2.5 for any
$a\in [a_1, a_2]$ in the sense that if $x_i$ is a solution of (\ref{2.2}) with Hamiltonian function $H_{a_i}(x)$ with
$i=1$ and $2$ respectively such that both $x_1$ and $x_2$ correspond to the same closed characteristic $(\tau,y)$ on
$\Sigma$, then we have
\be C_{S^1,\;q}(F_{a_1,K},\;S^1\cdot {x}_1) \cong C_{S^1,\;q}(F_{a_2,K},\;S^1\cdot{x}_2), \qquad \forall\,q\in \Z. \lb{3.4}\ee
In other words, the critical modules are independent of the choices of all $a>\frac{\tau}{T}$, the function $\vf_a$ satisfying
(i)-(ii) of Lemma 2.2, and $H_a$ satisfying Proposition 2.5. }

{\bf Proof.} Let $\vf_a$ be a family of functions satisfying (i)-(ii) of Lemma 2.2 and let $H_a(x)$ satisfy Proposition
2.5 parameterized by $a\in [a_1, a_2]$. Without loss of generality we can assume $H_a$ depends continuously on $a$ in the
sense of Remark 2.6. For each $a\in [a_1, a_2]$, we denote by $x_a$ the corresponding solution of (\ref{2.2}) with the
Hamiltonian $H_a$.

Now (\ref{2.26}) and (\ref{2.27}) imply that $b\mapsto\Psi_{b,K}$ is continuous in the $C^1$ topology. Then
$b\mapsto F_{b,K}$ is continuous in the $C^1$ topology too. Note that the number of critical orbits of each $F_{b,K}$
is finite. Hence by the continuity of critical modules (cf. Theorem 8.8 of \cite{MaW1} or Theorem 1.5.6 on p.53
of \cite{Cha1}, which can be easily generalized to the equivariant case), our proposition holds. Note that a similar
argument shows that the critical modules are independent of the choice of $\vf_a$ in $\wtd{H}_a(x)=a\vf_a(j(x))$
whenever $a$ is fixed, $\vf_a$ satisfies (i)-(ii) of Lemma 2.2, and $H_a$ satisfies Proposition 2.5. \hfill\hb

The rest of this section is devoted to the proof of independence of the critical modules for closed characteristics
in the choice of $K$. In the following, we fix an $a>\frac{\tau}{T}$, and write $F_K$ and $H$ for $F_{a, K}$ and $H_a$
respectively. We suppose also that $K\in {\bf R}$ satisfy (\ref{2.4}), i.e.,
\be  H_K(x) = H(x)+\frac{1}{2}K|x|^2 \qquad {\rm is\;strictly\;convex}. \lb{3.5}\ee

By Lemma 2.7, the critical points of $F_{K}$ which are the solutions of (\ref{2.2}) are the same for any
$K$ satisfying that $K\notin \frac{2\pi}{T}{\Z}$. Recall $d(K) = 2n([KT/{2\pi}]+1)$ by (\ref{2.15}).

{\bf Theorem 3.3.} {\it Suppose $\bar{x}$ is a nonzero critical point of $F_{K}$. Then the $S^1$-critical module
$C_{S^1,d(K)+l}(F_{K},S^1\cdot\bar{x})$ is independent of the choice of $K$ for $KT \notin 2\pi\Z$, i.e.,
\be C_{S^1,d(K)+l}(F_{K},S^1\cdot\bar{x}) \cong C_{S^1,d(K^\prime)+l}(F_{K^\prime},S^1\cdot\bar{x}),  \lb{3.6}\ee
where $KT$, $K^\prime T \notin 2\pi\Z$, $l\in \Z$, and both $K$ and $K^\prime$ satisfy (\ref{3.5}).}

We carry out the proof of this theorem in the following two cases: (a) $d(K)=d(K^\prime)$; and (b)
$K<K_0<K^\prime$ with $K_0-K, K^\prime-K_0$ small enough, and $\frac{K_0T}{2\pi}$ is an integer.

It is clear that proofs of Cases (a) and (b) imply the general case.

{\bf Proof of Case (a).} By Lemma 2.7, the critical orbits of $F_\sigma$ are independent of $\sigma$
and the number is finite by the condition (F). Then for $\bar{x}\in \crit(F_{\sg})$ there exists an
$S^1$-invariant open neighborhood $\mathcal{U}$ of $S^1\cdot\bar{x}$ in $X$ such that $F_{\sigma}$ has
a unique critical orbit $S^1\cdot\bar{x}$ in $\mathcal {U}$, for all $\sigma\in [K, K^\prime]$. By (7.12)
and (7.13) of \cite{Vit1}, $\sigma\mapsto F_{\sigma}$ is continuous in $C^1(\overline{\mathcal {U}})$
topology for $\sg\notin \frac{2\pi}{T}{\Z}$. Hence, by the continuity of critical modules (see Theorem 8.8
of \cite{MaW1} or Theorem 1.5.6 of \cite{Cha1}, which can be easily generalized to the equivariant sense),
and the Palais-Smale condition of $F_\sigma$ given by Proposition 2.12, the proof of Case (a) is complete.
\hfill\hb

Before we give the proof of Case (b), we give one definition and two lemmas.

{\bf Definition 3.4.} {\it Assume $K_0\in \frac{2\pi}{T}{\bf Z}$. As in \cite{Vit1}, let
$$  D_\infty(K_0)=\{\exp(-JK_0t)x\mid x\in {\bf R}^{2n}\}, $$
and define $C(K_0)$ to be the orthogonal complement of $D_\infty(K_0)$ in $X$.}

Let $K_0<K^\prime$ such that $K^\prime-K_0$ is small enough. As pointed out by Viterbo in \cite{Vit1},
$F_{K'}$ is strictly concave in the direction of $D_\infty(K_0)$. More precisely, similar to the argument to get
(7.7) of \cite{Vit1}, there exists a constant $C>0$ such that for all $x\in X$ and $h\in D_\infty(K_0)$ we have
\be \langle F_{K^\prime}^\prime(x+h)-F_{K^\prime}^\prime(x),h\rangle_X \le -C\|h\|_X^2,  \lb{3.7}\ee
where the equality holds if and only if $h=0$.

Let $\bar{x}$ be a nonzero critical point of $F_{K^\prime}$ with multiplicity $\mul(\bar{x})=m$, i.e., it corresponds
to a closed characteristic $(m\tau, y)\subset\Sigma$ with $(\tau, y)$ being prime. Write $\bar{x}=\bar{y}+\bar{z}$,
where $\bar{y}\in C(K_0)$, $\bar{z}\in D_\infty(K_0)$. Note that $\bar{y}\neq 0$ must hold, because otherwise, by
(7.14) of \cite{Vit1}, we have $\bar{z}=0$ and then $\bar{x}=0$.

Then expanding $\bar{x}$ into its Fourier series according to the $C(K_0)$ and $D_\infty(K_0)$ components, from the
fact $\mul(\bar{x})=m$, we obtain $\mul(\bar{y})=m$, then $\bar{y}(t+\frac{1}{m})=\bar{y}(t)$ for all $t\in \R$
and the orbit of $\bar{y}$, namely, $S^1\cdot \bar{y}\cong S^1/\Z_m\cong S^1$. Let
$p: N(S^1\cdot \bar{y})\rightarrow S^1\cdot \bar{y}$ be the normal bundle of $S^1\cdot \bar{y}$ in $C(K_0)$ and let
$p^{-1}(\theta\cdot \bar{y})=N(\theta\cdot \bar{y})$ be the fibre over $\theta\cdot \bar{y}$, where $\theta\in S^1$.
Let $DN(S^1\cdot \bar{y})$ be the $\varrho$ disk bundle of $N(S^1\cdot \bar{y})$ for some $\varrho>0$ sufficiently
small, i.e., $DN(S^1\cdot \bar{y})=\{\xi\in N(S^1\cdot \bar{y})\;| \; \|\xi\|_{H^1}<\varrho\}$ which is identified
by the exponential map with a subset of $C(K_0)$, and let
$DN(\theta\cdot \bar{y})=p^{-1}(\theta\cdot \bar{y})\cap DN(S^1\cdot \bar{y})$ be the disk over $\theta\cdot \bar{y}$.
Clearly, $DN(\theta\cdot \bar{y})$ is $\Z_m$-invariant and we have $DN(S^1\cdot \bar{y})=DN(\bar{y})\times_{\Z_m}S^1$,
where the $\Z_m$ action is given by
$$ (\theta, v, t)\in \Z_m\times DN(\bar{y})\times S^1\mapsto (\theta\cdot v, \;\theta^{-1}t)\in DN(\bar{y})\times S^1. $$
Hence for an $S^1$ invariant subset $\Gamma$ of $DN(S^1\cdot \bar{y})$, we have
$\Gamma/S^1=(\Gamma_{\bar{y}}\times_{\Z_m}S^1)/S^1=\Gamma_{\bar{y}}/\Z_m$, where $\Gamma_{\bar{y}}=\Gamma\cap DN(\bar{y})$.
Obviously, we also have a bundle $\tilde{p}: \wtd{N}(S^1\cdot \bar{x})\rightarrow S^1\cdot \bar{x}$ of $S^1\cdot \bar{x}$
in $X=C(K_0)\oplus D_\infty(K_0)$, where the fibre over $\theta\cdot \bar{x}$ is
$N(\theta\cdot \bar{y})\oplus D_\infty(K_0)$, $\theta\in S^1$. Let $D\wtd{N}(S^1\cdot \bar{x})$ be the $\varrho$ disk
bundle of $\wtd{N}(S^1\cdot \bar{x})$ for some $\varrho>0$ sufficiently small, i.e.,
$D\wtd{N}(S^1\cdot \bar{x})=\{\xi\in \wtd{N}(S^1\cdot \bar{x})\;| \; \|\xi\|_{H^1}<\varrho\}$ which is identified by the
exponential map with a subset of $X$, and let
$D\wtd{N}(\theta\cdot \bar{x})=\tilde{p}^{-1}(\theta\cdot \bar{x})\cap D\wtd{N}(S^1\cdot \bar{x})$ be the disk over
$\theta\cdot \bar{x}$. Clearly, $D\wtd{N}(\theta\cdot \bar{x})$ is $\Z_m$-invariant and we have
$D\wtd{N}(S^1\cdot \bar{x})=D\wtd{N}(\bar{x})\times_{\Z_m}S^1$ where the $\Z_m$ action is given by
$$ (\theta, v, t)\in \Z_m\times D\wtd{N}(\bar{x})\times S^1
         \mapsto (\theta\cdot v, \;\theta^{-1}t)\in D\wtd{N}(\bar{x})\times S^1. $$

{\bf Lemma 3.5.} {\it Let $\bar{x}$ be a nonzero critical point of $F_{K^\prime}$ with $\mul(\bar{x})=m$. Then there
exists an open ball $B(0,r)$ in $T_{\bar{x}}(D\wtd{N}(\bar{x}))$ centered at $0$ with radius
$r>0$, a local $\Z_m$-equivariant homeomorphism $\phi: B(0,r)\rightarrow \phi(B(0, r))\subset D\wtd{N}(\bar{x})$,
$\phi(0)=\bar{x}$, a $\Z_m$-equivariant $C^0$ map
$h: B(0, r)\cap T_{\bar{y}}(DN(\bar{y}))\rightarrow D_\infty(K_0)$ such that
\be F_{K^\prime}(\phi(\xi)) = -\|\eta\|_X^2+F_{K^\prime}(\bar{x}+\nu+h(\nu)), \lb{3.8} \ee
where $\xi=\eta+\nu$ with $\nu\in B(0,r)\cap T_{\bar{y}}(DN(\bar{y}))$ and $\eta\in D_\infty(K_0)$.}

{\bf Proof.} By (\ref{3.7}), in the direction of $D_\infty(K_0)$, $F_{K^\prime}$ is strictly concave.
Then there is a map $h: T_{\bar{y}}(DN(\bar{y}))\rightarrow D_\infty(K_0)$
uniquely defined by the relation $\nabla F_{K^\prime}(\bar{x}+\nu+h(\nu))\in C(K_0)$, i.e.,
$h(\nu)$ achieves the strict maximum of $F_{K^\prime}(\bar{x}+\nu+g)$ for $g$ in $D_\infty(K_0)$.
By the same proof of Lemma 2.2 of \cite{DHK1}, and noticing that (7.14) of \cite{Vit1},
we get that $h$ is continuous and $h(0)=0$. Note that since $F_{K^\prime}$ is $S^1$-invariant
and $\mul(\bar{x})=m$, then  $h$ is $\Z_m$-equivariant.

Let $H^+=T_{\bar{y}}(DN(\bar{y}))$ and $H^-=D_\infty(K_0)$. Now we define a map
$\psi:T_{\bar{x}}(D\wtd{N}(\bar{x}))=H^+\oplus H^-
\rightarrow T_{\bar{x}}(D\wtd{N}(\bar{x}))$ by\bea \psi(\nu, \mu)&=&(\nu, \psi_1(\nu, \mu))\nn\\&\equiv&
\left\{\matrix{(\nu, \sqrt{F_{K^\prime}(\bar{x}+\nu+h(\nu))-
F_{K^\prime}(\bar{x}+\nu+h(\nu)+\mu)}\frac{\mu}{\|\mu\|_X}), if~\mu\neq 0, \cr
(\nu, 0),~\qquad \qquad\qquad\qquad\qquad\qquad\qquad\qquad\qquad\qquad   if~\mu=0,\cr }\right.\lb{3.9} \eea
where $\nu\in H^+$, $\mu\in H^-$. Then $\psi$ is continuous on $H^+\oplus H^-$.

{\bf Claim (A).} {\it $\psi$ is one-to-one on $H^+\oplus H^-$.}

Suppose $\psi(\nu_1, \mu_1)=\psi(\nu_2, \mu_2)$ for some $(\nu_i, \mu_i)\in H^+\oplus H^-, i=1, 2$.
By (\ref{3.9}), we have $\nu_1=\nu_2$, $\frac{\mu_1}{\|\mu_1\|_X}=\frac{\mu_2}{\|\mu_2\|_X}$ and
\bea
&&\sqrt{F_{K^\prime}(\bar{x}+\nu_1+h(\nu_1))-F_{K^\prime}(\bar{x}+\nu_1+h(\nu_1)+\mu_1)}  \nn\\
&&\qquad = \sqrt{F_{K^\prime}(\bar{x}+\nu_2+h(\nu_2))-F_{K^\prime}(\bar{x}+\nu_2+h(\nu_2)+\mu_2)}. \nn\eea
Then we have $F_{K^\prime}(\bar{x}+\nu_1+h(\nu_1)+\mu_1)=F_{K^\prime}(\bar{x}+\nu_1+h(\nu_1)+\mu_2)$
and we may suppose $\mu_2=s \mu_1$ for some $s\geq 1$. By the mean value theorem, there exists
$1\leq t\leq s$ such that
\bea 0
&=& F_{K^\prime}(\bar{x}+\nu_1+h(\nu_1)+\mu_2)-F_{K^\prime}(\bar{x}+\nu_1+h(\nu_1)+\mu_1)\nn\\
&=& F_{K^\prime}(\bar{x}+\nu_1+h(\nu_1)+s\mu_1)-F_{K^\prime}(\bar{x}+\nu_1+h(\nu_1)+\mu_1)\nn\\
&=& \langle F_{K^\prime}^\prime(\bar{x}+\nu_1+h(\nu_1)+t\mu_1), (s-1)\mu_1\rangle_X\nn\\
&=& \langle F_{K^\prime}^\prime(\bar{x}+\nu_1+h(\nu_1)+t\mu_1)-F_{K^\prime}^\prime(\bar{x}+\nu_1+h(\nu_1)),
            (s-1)\mu_1\rangle_X\nn\\&\leq&-(s-1) tC\|\mu_1\|_X^2£¬ \nn\eea
where $C>0$ is defined in (\ref{3.7}), and we have used the fact that
$\nabla F_{K^\prime}(\bar{x}+\nu_1+h(\nu_1))\in C(K_0)$ implies
$$ \langle F_{K^\prime}^\prime(\bar{x}+\nu_1+h(\nu_1)), \mu_1\rangle_X = 0. $$
Thus we have $s=1$ or $\mu_1=0$, then $\mu_1=\mu_2$. Claim (A) follows.

{\bf Claim (B).} {\it For any $\epsilon>0$, there exists a positive real number $\delta_\epsilon>0$ such that
\bea B_{H^+}(0, \epsilon)\times B_{H^-}(0, \delta_\epsilon)\subseteq
\psi(B_{H^+}(0, \epsilon)\times B_{H^-}(0, \epsilon)),\nn\eea
where $B_{H^*}(0, \epsilon)$ denotes an open ball in $H^*$, $*=+, -$.}

In fact, by (\ref{3.7}) and noticing that $\nabla F_{K^\prime}(\bar{x}+\nu+h(\nu))\in C(K_0)$ implies
$$ \langle F_{K^\prime}^\prime(\bar{x}+\nu+h(\nu)), \mu\rangle_X=0, \qquad \forall \mu\in D_\infty(K_0), $$
then we have
\bea
&&F_{K^\prime}(\bar{x}+\nu+h(\nu))-F_{K^\prime}(\bar{x}+\nu+h(\nu)+\mu)  \nn\\
&&\quad = -\int_0^1 \frac{d}{dt}F_{K^\prime}(\bar{x}+\nu+h(\nu)+t\mu) dt  \nn\\
&&\quad = -\int_0^1 \langle F_{K^\prime}^\prime(\bar{x}+\nu+h(\nu)+t\mu), \mu\rangle_X dt  \nn\\
&&\quad = -\int_0^1 \langle F_{K^\prime}^\prime(\bar{x}+\nu+h(\nu)+t\mu)
               - F_{K^\prime}^\prime(\bar{x}+\nu+h(\nu)), \mu\rangle_X dt \nn\\
&&\quad \ge C\int_0^1 t\|\mu\|_X^2 dt=\frac{C}{2}\|\mu\|_X^2,  \nn\eea
where $C>0$ is the constant defined in (\ref{3.7}). Then by the continuity of $h$ and the definition of
$\psi_1$ in (\ref{3.9}), it follows that
$$ \{t\frac{\mu}{\|\mu\|_X}, 0\leq t < \sqrt{\frac{C}{2}}\epsilon\}
  \subseteq \psi_1(\nu\times B_{H^-}(0, \epsilon)), \quad
      \forall \mu\in \partial B_{H^-}(0, \epsilon), \nu\in B_{H^+}(0, \epsilon). $$
Let $\delta_\epsilon=\sqrt{\frac{C}{2}}\epsilon$, we have
$$  B_{H^+}(0, \epsilon)\times B_{H^-}(0, \delta_\epsilon)
      \subseteq \psi(B_{H^+}(0, \epsilon)\times B_{H^-}(0, \epsilon)), $$
which proves Claim (B).

Let $\varphi$ be the restriction of $\psi^{-1}$ on $B_{H^+}(0, \epsilon)\times B_{H^-}(0, \delta_\epsilon)$.

{\bf Claim (C).} {\it $\varphi$ is continuous on $B_{H^+}(0, \epsilon)\times B_{H^-}(0, \delta_\epsilon)$.}

Let $(\nu_0, \mu_0)\in \psi^{-1}(B_{H^+}(0, \epsilon)\times B_{H^-}(0, \delta_\epsilon))$ and
$\{(\nu_n, \mu_n)\}$ be a sequence in $\psi^{-1}(B_{H^+}(0, \epsilon)\times B_{H^-}(0, \delta_\epsilon))$ such that
$\{\psi(\nu_n, \mu_n)\}$ converges to $\psi(\nu_0, \mu_0)$, we prove that $\{(\nu_n, \mu_n)\}$ converges
to $(\nu_0, \mu_0)$. Firstly, by definition we have $\nu_n\rightarrow \nu_0$.
Since $\overline{B_{H^-}(0, \epsilon)}$ is compact, we can suppose
that $\mu_n\rightarrow \mu$ in $\overline{B_{H^-}(0, \epsilon)}$. Then by the continuity of $h$ and $F_{K^\prime}$,
$h(\nu_n)\rightarrow h(\nu_0)$ and
\bea
&&\sqrt{F_{K^\prime}(\bar{x}+\nu_n+h(\nu_n))-F_{K^\prime}(\bar{x}+\nu_n+h(\nu_n)+\mu_n)}  \nn\\
&&\qquad\rightarrow \sqrt{F_{K^\prime}(\bar{x}+\nu_0+h(\nu_0))-F_{K^\prime}(\bar{x}+\nu_0+h(\nu_0)+\mu)}
           \qquad {\rm as} \quad n\to +\infty. \nn\eea
Thus we get $\psi(\nu_0, \mu_0)=\psi(\nu_0, \mu)$ by (\ref{3.9}), since $\psi$ is one-to-one,
then $\mu=\mu_0$ and Claim (C) follows.

Now by Claims (A), (B) and (C), $\varphi$ is an homeomorphism from
$B_{H^+}(0, \epsilon)\times B_{H^-}(0, \delta_\epsilon)$ to
an open neighborhood of 0 in $T_{\bar{x}}(D\wtd{N}(\bar{x}))$. We define $\phi=\exp_{\bar{x}}\circ\tau\circ\varphi$,
where $\exp_{\bar{x}}$ is the exponential map, $\tau$ is defined by
\bea \tau(\nu, \mu)=(\nu, h(\nu)+\mu),\nn\eea
which is an homeomorphism from $T_{\bar{x}}(D\wtd{N}(\bar{x}))$ to itself satisfying $\tau(0)=0$.
Then for $\epsilon>0$ sufficiently small,
$\phi$ is an homeomorphism from $B_{H^+}(0, \epsilon)\times B_{H^-}(0, \delta_\epsilon)$ to
an open neighborhood of $\bar{x}$ in $D\wtd{N}(\bar{x})$. Note that by the above proof, $\phi$ is $\Z_m$-equivariant,
and for any $(\nu, \eta)\in B_{H^+}(0, \epsilon)\times B_{H^-}(0, \delta_\epsilon)$, we can write
$\eta=\psi_1(\nu, \mu)$ for some $\mu\in B_{H^-}(0, \epsilon)$, then
\bea F_{K^\prime}(\phi(\nu, \eta))
&=& F_{K^\prime}(\phi\circ \psi(\nu, \mu)) = F_{K^\prime}(\bar{x}+\nu+h(\nu)+\mu)  \nn\\
&=& F_{K^\prime}(\bar{x}+\nu+h(\nu))-\|\psi_1(\nu, \mu)\|_X^2 = F_{K^\prime}(\bar{x}+\nu+h(\nu))-\|\eta\|_X^2.  \nn\eea
Then (\ref{3.8}) holds. Let $r=\min\{\epsilon, \delta_\epsilon\}$, we complete the proof of Lemma 3.5. \hfill\hb

{\bf Lemma 3.6.} {\it $F_{K_0}$ satisfies the Palais-Smale condition on $C(K_0)$.}

{\bf Proof.} Similar to the study in Section 3 of \cite{Vit1}, since the map $y\to -J\dot{y}+K_0 y$ is a Hilbert
space isomorphism between $C(K_0)$ and a subspace $E_0$ of $E=L^2(\R/(T\Z),{\bf R}^{2n})$, we can define a functional
$\Psi_{K_0}$ on $E_0$ as in (\ref{2.9}) (cf. also (3.5) of \cite{Vit1}). Then the corresponding Propositions 3.9 and
4.1 of \cite{Vit1} hold in our case. Here the equation $M(g_n+h_n)-\nabla H^*_K(g_n+h_n)=\ep_n$ in the proof of
Proposition 4.1 of \cite{Vit1} should be modified to $M(g_n+h_n)-\nabla H^*_K(g_n+h_n)=\ep_n+e_n$ with
$e_n\in D_\infty(K_0)$ and $\ep_n\to 0$, where $M$ is the inverse of the map $y\to -J\dot{y}+K_0 y$ on $C(K_0)$,
$z_n=Mg_n+Mh(g_n)$ should be modified to $z_n=Mg_n+Mh(g_n)-e_n$. As in the proof of Proposition 4.1 of \cite{Vit1},
this $z_n$ also tends to infinity in the $C^0$ topology. Then by the proof of Proposition 4.1 of \cite{Vit1} we obtain
the Palais-Smale condition for $\Psi_{K_0}$. Now using the same argument of our Proposition 2.12, we obtain that
$F_{K_0}$ satisfies the Palais-Smale condition on $C(K_0)$. \hfill\hb

We now continue our proof of Theorem 3.3.

{\bf Proof of Case (b).} Similar to our discussion in (\ref{3.7}) and Lemma 3.5 for the special case $L=K'$, by the
proof of (7.6) and (7.7) in \cite{Vit1}, the functional $F_L$ is strictly convex (resp. concave) in the direction of
$D_\infty(K_0)$ for $L<K_0$ (resp. $L>K_0$). Thus there is an $S^1$-equivariant map
\be  z_L: C(K_0)\to D_\infty(K_0), \qquad y\mapsto z_L(y),  \lb{3.10}\ee
uniquely defined by the relation $\nabla F_L(y+z_L(y))\in C(K_0)$, i.e., $z_L(y)$ achieves the minimum (resp. maximum)
of $F_L(y+h)$ for $h$ in $D_\infty(K_0)$. Note that the map $z_L$ is $C^{0,1}$ by Page 627 of \cite{Vit1}, and when
$L=K'$ the map $z_L$ is the map $h$ defined in Lemma 3.5.

Now for $y\in C(K_0)$, set $\wtd{F}_L(y)=F_L(y+z_L(y))$ for $L\in [K,K']\bs\{K_0\}$. Then $\wtd{F}_L$ is $C^{1,1}$ and
$S^1$-invariant. By Lines 9-10 on Page 642 of \cite{Vit1}, we have
\bea
|\wtd{F}_L(y)-\wtd{F}_{K_0}(y)| &\le& C|L-K_0|\|y\|_{H^1}^2,   \lb{3.11}\\
\|\nabla\wtd{F}_L(y)-\nabla\wtd{F}_{K_0}(y)\|_{H^1} &\le& C|L-K_0|\|y\|_{H^1}^2,   \lb{3.12}\eea
where we write $\wtd{F}_{K_0}$ for $F_{K_0}$ restricted to $C(K_0)$. Since $\bar{x}\neq 0$ is a critical point of $F_L$,
there exists a unique $\bar{y}\in C(K_0)$ such that $\bar{x}=\bar{y}+\bar{z}$, where $\bar{z}\in D_\infty(K_0)$. Then
$z_L(\bar{y})=\bar{z}$ holds, which is independent of $L$. Note that $\bar{y}\neq 0$, because otherwise, by (7.14) of
\cite{Vit1}, we have $\bar{z}=0$ and then $\bar{x}=0$.

To continue the proof, we need the following three claims.

{\bf Claim 1.} {\it We can choose an $S^1$-invariant open neighborhood $\mathcal{U}_1$ of $S^1\cdot\bar{y}$ in
$C(K_0)$ such that $S^1\cdot\bar{y}$ is the unique critical orbit of $\wtd{F}_L$ in $\mathcal {U}_1$ for all
$L\in [K, K^\prime]$.}

In fact, firstly there is an $S^1$-invariant open neighborhood $V_1$ of $S^1\cdot\bar{y}$ in $C(K_0)$ such that
$S^1\cdot \bar{y}$ is the unique critical orbit of $\wtd{F}_L$ in $V_1$ for all $L\in [K, K^\prime]\bs K_0$.
If not, there exists a sequence of $S^1$ orbits $\{S^1\cdot \bar{y}_j\}_{j\geq 1}\subset C(K_0)$ such that
$\lim_{j\to\infty}{\bar{y}_j}=\bar{y}$ and $\bar{y}_j$ is a critical point of $\wtd{F}_{L_j}$ for some
$L_j\in [K, K^\prime]\bs K_0$. Then $\bar{y}_j+z_{L_j}(\bar{y}_j)$ is a critical point of $F_{L_j}$ for all
$j\in {\bf N}$. Since by Lemma 2.7 the critical points of $F_{L}$ are the same for all
$L\in [K, K^\prime]\setminus K_0$, then $\bar{y}_j+z_{L_j}(\bar{y}_j)$ is a critical point of $F_L$ for some
$L\in [K,K']\bs\{K_0\}$ and $z_{L_j}(\bar{y}_j)=z_{L}(\bar{y}_j)$. But
$\lim_{j\to\infty}{(\bar{y}_j+z_L(\bar{y}_j))}=\bar{y}+z_L(\bar{y})$ and the critical orbits of $F_L$ are isolated
by the condition (F), which yields a contradiction.

If $\nabla\wtd{F}_{K_0}(z)=0$ for some $z\in C(K_0)$, by definition, we have
$w=z-\nabla H_{K_0}^*(-J\dot{z}+K_0 z)\in D_\infty(K_0)$. Then
$\nabla H_{K_0}(z-w)=-J\dot{z}+K_0 z=-J(\dot{z}-\dot{w})+K_0(z-w)$, i.e., $\nabla H(z-w)=-J(\dot{z}-\dot{w})$.
Thus $z-w$ is a solution of (2.2). But on the other hand, all the solutions of (2.2) are isolated $S^1$ orbits by
the condition (F), so we can choose an $S^1$-invariant open neighborhood $V_2$ of $\bar{y}$ in $C(K_0)$ such that
$S^1\cdot \bar{y}$ is the unique critical orbit of $\wtd{F}_{K_0}$ in $V_2$. Hence, setting
$\mathcal{U}_1=V_1\cap V_2$, Claim 1 is proved.

Note that $\wtd{F}_L$ satisfies the Palais-Smale condition by Proposition 2.12 and Lemma 3.6. Now combining
(\ref{3.11})-(\ref{3.12}) with the continuity of critical modules depending on $L$ (cf. Theorem 8.8 of \cite{MaW1}
or Theorem 1.5.6 of \cite{Cha1}, which can be easily generalized to the equivariant sense), we obtain the
$C^1$-continuity of $\td{F}_L$ in $L\in [K,K']$. Together with Claim 1, we obtain
\be  C_{S^1,d(K)+l}(\wtd{F}_{K},S^1\cdot\bar{y}) \cong C_{S^1,d(K)+l}(\wtd{F}_{K_0},S^1\cdot\bar{y})
     \cong C_{S^1,d(K)+l}(\wtd{F}_{K^\prime},S^1\cdot\bar{y}).   \lb{3.13}\ee

{\bf Claim 2.} {\it $F_{L}(\bar{x})$ is independent of $L$.}

In fact, since $\nabla H_L(\bar{x})=-J\dot{\bar{x}}+L\bar{x}$, then
$H_L^*(-J\dot{\bar{x}}+L\bar{x})=(-J\dot{\bar{x}}+L\bar{x}, \bar{x})-H_L(\bar{x})$. Thus
\bea F_{L}(\bar{x})
&=& \int_0^T{\left[\frac{1}{2}(J\dot{\bar{x}}-L\bar{x},\bar{x})+H_L^*(-J\dot{\bar{x}}+L\bar{x})\right]dt}  \nn\\
&=& \int_0^T{\left[-\frac{1}{2}(J\dot{\bar{x}}-L\bar{x},\bar{x})-H_L(\bar{x})\right]dt}  \nn\\
&=& \int_0^T{\left[-\frac{1}{2}(J\dot{\bar{x}},\bar{x})-H(\bar{x})\right]dt},  \lb{3.14}\eea
which is independent of $L$. Thus Claim 2 is proved.

Now let $c=F_{L}(\bar{x})$. We define $X^c(K)=\{y\in X\mid F_K(y)\leq c\}$ and
$\wtd{X}^c(K)=\{y\in C(K_0)\mid \wtd{F}_K(y)\leq c\}$. Let $\wtd{U}$ be an $S^1$-invariant open neighborhood of
$\bar{y}$ in $C(K_0)$ such that $\wtd{F}_K$ has unique critical orbit $S^1\cdot \bar{y}$ in $\wtd{U}$, then
$U\equiv\wtd{U}\times D_\infty(K_0)$ is an $S^1$-invariant open neighborhood of $\bar{x}$ such that $F_K$ has unique
critical orbit $S^1\cdot \bar{x}$ in $U$.

{\bf Claim 3.} {\it The natural embeddings $\wtd{X}^c(K)\cap\wtd{U} \to X^c(K)\cap U$ and
$(\wtd{X}^{c}(K)\bs\{\bar{y}\})\cap\wtd{U} \to (X^{c}(K)\bs\{\bar{x}\})\cap U$ are $S^1$-equivariant homotopy
equivalences.}

In fact, by the strictly convexity of $F_K$ in the direction of $D_\infty(K_0)$ and the argument of Lemma
5.1 of \cite{Vit1}, the claim follows.

By Claim 3, we have the following:
\be  C_{S^1, d(K)+l}(F_{K}, S^1\cdot\bar{x})  \cong C_{S^1, d(K)+l}(\wtd{F}_{K}, S^1\cdot\bar{y}). \lb{3.15}\ee
Together with (\ref{3.13}), this yields
\be C_{S^1, d(K)+l}(F_{K}, S^1\cdot\bar{x}) \cong C_{S^1, d(K)+l}(\wtd{F}_{ K^\prime}, S^1\cdot\bar{y}).  \lb{3.16}\ee

In Lemma 3.5, let $f_1(\eta)=-\|\eta\|_X^2$ for all $\eta\in D_\infty(K_0)$,
$f_2(\nu)=F_{K^\prime}(\bar{x}+\nu+h(\nu))$ for all $\nu\in B(0,r)\cap T_{\bar{y}}(DN(\bar{y}))$.
Then the Gromoll-Meyer pair of 0 for $f_1$, i.e., $(W_1, W_{1^-})$, is
$\Z_m$-equivariant homotopy equivalent with $(B^{2n}, S^{2n-1})$ since $f_1$ is $\Z_m$-invariant. Note that for (\ref{3.8})
of Lemma 3.5, by the definitions of $h$ and $\wtd{F}_{K^\prime}$, we have
\bea
&& \frac{\partial}{\partial z}F_{K^\prime}(\bar{y}+\nu+z_{K^\prime}(\bar{y}+\nu))=0, \nn\\
&& f_2(\nu) = F_{K^\prime}(\bar{x}+\nu+h(\nu))
         = F_{K^\prime}(\bar{y}+\nu+z_{K^\prime}(\bar{y}+\nu))=\wtd{F}_{K^\prime}(\bar{y}+\nu)£¬ \nn\eea
for $\nu\in T_{\bar{y}}(DN(\bar{y}))$. Denote
by $(W_2, W_{2^-})$ the Gromoll-Meyer pair of $\bar{y}$ with respect to the negative gradient vector field of $\wtd{F}_{K^\prime}$
in $DN(\bar{y})$, $(W_2, W_{2^-})$ is $\Z_m$-invariant since $f_2$ is $\Z_m$-invariant. Thus, we obtain
\be  C_{d(K)+l}(\wtd{F}_{K^\prime}|_{DN(\bar{y})}, \bar{y})
       \cong H_{d(K)+l}(W_2, W_{2^-}).  \lb{3.17}\ee
Using Lemma 1.5.1 of \cite{Cha1} and Lemma 3.5, we have
\be C_{d(K)+l+2n}(F_{K^\prime}|_{D\wtd{N}(\bar{x})}, \bar{x})
   \cong H_{d(K)+l+2n}(W_1\times W_2, (W_{1^-}\times W_2)\cup(W_1\times W_{2^-})). \lb{3.18}\ee
By Definition 3.1, we have
\bea
&& C_{S^1,\;\ast}(F_{K^\prime}, \;S^1\cdot \bar{x}) \nn\\
&&\qquad\qquad \cong
   H_{S^1, \;\ast}(X_{a,K^\prime}(\bar{x})\cap D\wtd{N}(S^1\cdot\bar{x}),
      \;(X_{a, K^\prime}(\bar{x})\setminus (S^1\cdot \bar{x}))\cap D\wtd{N}(S^1\cdot\bar{x})),  \lb{3.19}\eea
where $X_{a, K^\prime}(\bar{x})$ is defined as in (\ref{3.2}). Since all the isotropy groups
$G_x=\{g\in S^1\;|\;g\cdot x=x\}$ for $x\in D\wtd{N}(S^1\cdot \bar{x})$ are finite, we can use Lemma 6.11 of
\cite{FaR1} to obtain
\bea
&& H_{S^1}^\ast(X_{a, K^\prime}(\bar{x})\cap D\wtd{N}(S^1\cdot\bar{x}),
      \; (X_{a,K^\prime}(\bar{x})\bs (S^1\cdot \bar{x}))\cap D\wtd{N}(S^1\cdot\bar{x}))  \nn\\
&&\qquad \cong H^\ast (X_{a, K^\prime}(\bar{x})\cap D\wtd{N}(S^1\cdot\bar{x})/S^1,\;
      (X_{a,K^\prime}(\bar{x})\setminus (S^1\cdot \bar{x}))\cap D\wtd{N}(S^1\cdot\bar{x})/S^1) \nn\\
&&\qquad \cong H^\ast (X_{a, K^\prime}(\bar{x})\cap D\wtd{N}(\bar{x})/\Z_m,
      \; (X_{a,K^\prime}(\bar{x})\bs(\bar{x}))\cap D\wtd{N}(\bar{x})/\Z_m). \nn\eea
By the condition (F) at the beginning of Section 2, a small perturbation on the energy functional
can be applied to reduce each critical orbit to nearby non-degenerate ones. Thus similar to the
proofs of Lemma 2 of \cite{GrM1} and Lemma 4 of \cite{GrM2}, all the homological $\Q$-modules of each space
pair in the above relations are all finitely generated. Therefore we can apply Theorem 5.5.3 and
Corollary 5.5.4 on pages 243-244 of \cite{Spa1} to obtain the same relation on homological $\Q$-modules:
\bea
&& H_{S^1,\ast}(X_{a, K^\prime}(\bar{x})\cap D\wtd{N}(S^1\cdot\bar{x}),
      \; (X_{a,K^\prime}(\bar{x})\setminus (S^1\cdot \bar{x}))\cap D\wtd{N}(S^1\cdot\bar{x})) \nn\\
&&\qquad \cong H_\ast (X_{a, K^\prime}(\bar{x})\cap D\wtd{N}(S^1\cdot\bar{x})/S^1,
      \; (X_{a,K^\prime}(\bar{x})\setminus (S^1\cdot \bar{x}))\cap D\wtd{N}(S^1\cdot\bar{x})/S^1) \nn\\
&&\qquad \cong H_\ast (X_{a, K^\prime}(\bar{x})\cap D\wtd{N}(\bar{x})/\Z_m,
      \; (X_{a,K^\prime}(\bar{x})\setminus (\bar{x}))\cap D\wtd{N}(\bar{x})/\Z_m).  \lb{3.20}\eea
For a $\Z_m$-space pair $(A, B)$, let
$$ H_{\ast}(A, B)^{\pm\Z_m} = \{\sigma\in H_{\ast}(A,B)\,|\,L_{\ast}\sigma=\pm \sigma\}, $$
where $L$ is a generator of the $\Z_m$-action. Note that the same argument as in Section 6.3 of \cite{Rad2},
in particular Satz 6.6 of \cite{Rad2}, Lemma 3.6 of \cite{BaL1} or Theorem 3.2.4 of \cite{Bre1} yields
\bea
&& H_\ast(X_{a, K^\prime}(\bar{x})\cap D\wtd{N}(\bar{x})/\Z_m,
      \; (X_{a,K^\prime}(\bar{x})\setminus (\bar{x}))\cap D\wtd{N}(\bar{x})/\Z_m)  \nn\\&&
\qquad \cong H_\ast(X_{a, K^\prime}(\bar{x})\cap D\wtd{N}(\bar{x}),
      \;(X_{a, K^\prime}(\bar{x})\setminus (\bar{x}))\cap D\wtd{N}(\bar{x}))^{\Z_m}. \lb{3.21}\eea
Combining (\ref{3.19})-(\ref{3.21}), we have
\bea
&& C_{S^1,\;d(K)+l+2n}(F_{K^\prime}, \;S^1\cdot \bar{x}) \nn\\
&&\qquad \cong H_{d(K)+l+2n} (X_{a,K^\prime}(\bar{x})\cap D\wtd{N}(\bar{x}),
      \; (X_{a,K^\prime}(\bar{x})\setminus (\bar{x}))\cap D\wtd{N}(\bar{x}))^{\Z_m}  \nn\\
&&\qquad \cong C_{d(K)+l+2n}(F_{K^\prime}|_{D\wtd{N}(\bar{x})}, \bar{x})^{\Z_m}.   \lb{3.22}\eea
Similarly, we have
\be C_{S^1, d(K)+l}(\wtd{F}_{K^\prime},S^1\cdot\bar{y})
     \cong C_{d(K)+l}(\wtd{F}_{K^\prime}|_{DN(\bar{y})}, \bar{y})^{\Z_m}. \lb{3.23}\ee
Now by (\ref{3.17}) and (\ref{3.18}), as in Proposition 3.10 of \cite{WHL1}, we have
\be  C_{d(K)+l+2n}(F_{K^\prime}|_{D\wtd{N}(\bar{x})}, \bar{x})^{\Z_m}\cong
      C_{d(K)+l}(\wtd{F}_{K^\prime}|_{DN(\bar{y})}, \bar{y})^{\Z_m}.  \lb{3.24}\ee
In fact, let $\theta$ be a generator of the linearized $\Z_m$-action on $W_1$. Then $\theta(\xi)=\xi$ for
$0\neq \xi \in T_0(W_1)$ if and only if $m|\frac{K_0T}{2\pi}$. Thus together with (\ref{3.17}), (\ref{3.18})
and the fact that $\dim W_1$ is even, it yields (\ref{3.24}).

Hence, it follows from (\ref{3.22})-(\ref{3.24}) that
\be C_{S^1,d(K)+l+2n}(F_{K^\prime}, S^1\cdot\bar{x})  \cong  C_{S^1,d(K)+l}(\wtd{F}_{K^\prime}, S^1\cdot\bar{y}). \lb{3.25}\ee

Combining (\ref{3.16}) and (\ref{3.25}), using the fact that $d(K^\prime)=d(K)+2n$, we obtain
$$ C_{S^1, d(K)+l}(F_{K}, S^1\cdot\bar{x}) \cong C_{S^1, d(K)+l+2n}(F_{K^\prime}, S^1\cdot\bar{x})
      = C_{S^1,d(K^{\prime})+l}(F_{K^\prime}, S^1\cdot\bar{x}).  $$
The proof of Theorem 3.3 is complete. \hfill\hb

\setcounter{equation}{0}
\section{Periodic property of critical modules for closed characteristics}

In this section, we fix $a$ and let $u_K\neq 0$ be a critical point of $\Psi_{a, K}$ with multiplicity $\mul(u_K)=m$,
that is, $u_K$ corresponds to a closed characteristic $(\tau, y)\subset\Sigma$ with $(\tau, y)$ being $m$-iteration of
some prime closed characteristic. Precisely, by Proposition 2.5 and Lemma 2.7, we have $u_K=-J\dot x+Kx$ with $x$
being a solution of (\ref{2.2}) and $x=\rho y(\frac{\tau t}{T})$ with $\frac{\vf_a^\prime(\rho)}{\rho}=\frac{\tau}{aT}$.
Moreover, $(\tau, y)$ is a closed characteristic on $\Sigma$ with minimal period $\frac{\tau}{m}$. Hence the isotropy
group satisfies $\{\theta\in S^1\;|\;\theta\cdot u_K=u_K\}=\Z_m$ and the orbit of $u_K$, namely,
$S^1\cdot u_K\cong S^1/\Z_m\cong S^1$. By Lemma 2.10, we obtain a critical point $g_K$ of $\psi_{a,K}$ corresponding to
$u_K$, and then the isotropy group satisfies $\{\theta\in S^1\;|\;\theta\cdot g_K=g_K\}=\Z_m$.
Let $p: N(S^1\cdot g_K)\rightarrow S^1\cdot g_K$ be the normal bundle of $S^1\cdot g_K$ in $G$(as defined in Lemma 2.10)
and let $p^{-1}(\theta\cdot g_K)=N(\theta\cdot g_K)$ be the fibre over $\theta\cdot g_K$, where $\theta\in S^1$. Let
$DN(S^1\cdot g_K)$ be the $\varrho$ disk bundle of $N(S^1\cdot g_K)$ for some $\varrho>0$ sufficiently small, i.e.,
$DN(S^1\cdot g_K)=\{\xi\in N(S^1\cdot g_K)\;| \; \|\xi\|<\varrho\}$ which is identified by the exponential map with a
subset of $G$, and let $DN(\theta\cdot g_K)=p^{-1}(\theta\cdot g_K)\cap DN(S^1\cdot g_K)$ be the disk over
$\theta\cdot g_K$. Clearly, $DN(\theta\cdot g_K)$ is $\Z_m$-invariant and we have
$DN(S^1\cdot g_K)=DN(g_K)\times_{\Z_m}S^1$ where the $Z_m$ action is given by
$$(\theta, v, t)\in \Z_m\times DN(g_K)\times S^1\mapsto (\theta\cdot v, \;\theta^{-1}t)\in DN(g_K)\times S^1. $$
Hence for an $S^1$ invariant subset $\Gamma$ of $DN(S^1\cdot g_K)$, we have
$\Gamma/S^1=(\Gamma_{g_K}\times_{\Z_m}S^1)/S^1=\Gamma_{g_K}/\Z_m$, where $\Gamma_{g_K}=\Gamma\cap DN(g_K)$.

For a $\Z_m$-space pair $(A, B)$, let
\be H_{\ast}(A, B)^{\pm\Z_m} = \{\sigma\in H_{\ast}(A, B)\,|\,L_{\ast}\sigma=\pm \sigma\}, \lb{4.1}\ee
where $L$ is a generator of the $\Z_m$-action. Then as in Section 6 of \cite{Rad2}, Section 3 of \cite{BaL1} or
Lemma 3.9 of \cite{WHL1}, we have

{\bf Lemma 4.1.} {\it Suppose $u_K\neq 0$ is a critical point of $\Psi_{a, K}$ with $\mul(u_K)=m$, $g_K$ is a critical
point of $\psi_{a,K}$ corresponding to $u_K$. Then we have}
\bea C_{S^1,\;\ast}(\Psi_{a,K}, \;S^1\cdot u_K)
&\cong&C_{S^1,\; \ast}(\psi_{a,K},\;S^1\cdot g_K) \nn\\
&\cong& H_\ast((\wtd{\Lambda}_{a,K}(g_K)\cap DN(g_K))/\Z_m,\;
    ((\wtd{\Lambda}_{a,K}(g_K)\setminus\{g_K\})\cap DN(g_K))/\Z_m) \nn\\
&\cong& H_\ast((\wtd{\Lambda}_{a,K}(g_K)\cap DN(g_K)),\;
    ((\wtd{\Lambda}_{a,K}(g_K)\setminus\{g_K\})\cap DN(g_K)))^{\Z_m}.   \lb{4.2}\eea
where $\wtd{\Lambda}_{a,K}(g_K)=\{g\in G \;|\; \psi_{a, K}(g)\le\psi_{a, K}(g_K)\}$.

{\bf Proof.} For reader's conveniences, we sketch a proof here and refer to Section 6 of \cite{Rad2},
Section 3 of \cite{BaL1} or Lemma 3.9 of \cite{WHL1} for related details.

By Lemma 2.10 (iv), we have
$$ C_{S^1,\; \ast}(\Psi_{a,K},\;S^1\cdot u_K) \cong C_{S^1,\; \ast}(\psi_{a,K}, \;S^1\cdot g_K). $$
By Definition 3.1, we have
$$ C_{S^1,\;\ast}(\psi_{a,K}, \;S^1\cdot g_K)
   \cong H_{S^1,\;\ast}(\wtd{\Lambda}_{a,K}(g_K)\cap DN(S^1\cdot g_K),
       \;(\wtd{\Lambda}_{a,K}(g_K)\setminus (S^1\cdot g_K))\cap DN(S^1\cdot g_K)). $$
Since all the isotropy groups $A_x=\{a\in S^1\;|\;a\cdot x=x\}$ for $x\in DN(S^1\cdot g_K)$ are finite, we can use
Lemma 6.11 of \cite{FaR1} to obtain
\bea
&& H_{S^1}^\ast(\wtd{\Lambda}_{a,K}(g_K)\cap DN(S^1\cdot g_K),\;
          (\wtd{\Lambda}_{a,K}(g_K)\setminus (S^1\cdot g_K))\cap DN(S^1\cdot g_K)) \nn\\
&&\qquad \cong H^\ast ((\wtd{\Lambda}_{a,K}(g_K)\cap DN(S^1\cdot g_K))/S^1, \;
          ((\wtd{\Lambda}_{a,K}(g_K)\setminus (S^1\cdot g_K))\cap DN(S^1\cdot g_K))/S^1) \nn\\
&&\qquad \cong  H^\ast((\wtd{\Lambda}_{a,K}(g_K)\cap DN(g_K))/\Z_m,\;
    ((\wtd{\Lambda}_{a,K}(g_K)\setminus\{g_K\})\cap DN(g_K))/\Z_m).  \nn\eea
By the condition (F) at the beginning of Section 2, a small perturbation on the energy functional can be applied
to reduce each critical orbit to nearby non-degenerate ones. Thus similar to the proofs of Lemmas 2 and 4 of
\cite{GrM1}, all the homological $\Q$-modules of each space pair in the above relations are all finitely generated.
Therefore we can apply Theorem 5.5.3 and Corollary 5.5.4 on pages 243-244 of \cite{Spa1} to obtain the same relation
on homological $\Q$-modules:
\bea
&& H_{S^1,*}(\wtd{\Lambda}_{a,K}(g_K)\cap DN(S^1\cdot g_K),\;
    (\wtd{\Lambda}_{a,K}(g_K)\setminus (S^1\cdot g_K))\cap DN(S^1\cdot g_K))  \nn\\
&&\qquad \cong H_*((\wtd{\Lambda}_{a,K}(g_K)\cap DN(S^1\cdot g_K))/S^1, \;
    ((\wtd{\Lambda}_{a,K}(g_K)\setminus (S^1\cdot g_K))\cap DN(S^1\cdot g_K))/S^1)  \nn\\
&&\qquad \cong  H_\ast((\wtd{\Lambda}_{a,K}(g_K)\cap DN(g_K))/\Z_m,\;
    ((\wtd{\Lambda}_{a,K}(g_K)\setminus\{g_K\})\cap DN(g_K))/\Z_m).  \nn\eea
Note that the same argument as in Section 6.3 of \cite{Rad2}, in particular Satz 6.6 of \cite{Rad2}, Lemma
3.6 of \cite{BaL1} or Theorem 3.2.4 of \cite{Bre1} yields
\bea
&& H_\ast((\wtd{\Lambda}_{a,K}(g_K)\cap DN(g_K))/\Z_m,\;
    ((\wtd{\Lambda}_{a,K}(g_K)\setminus\{g_K\})\cap DN(g_K))/\Z_m) \nn\\
&&\qquad \cong H_\ast((\wtd{\Lambda}_{a,K}(g_K)\cap DN(g_K)),\;
    ((\wtd{\Lambda}_{a,K}(g_K)\setminus\{g_K\})\cap DN(g_K)))^{\Z_m}. \nn\eea
The above relations together complete the proof of Lemma 4.1. \hfill\hb

By (\ref{2.6}) and (\ref{2.9}), we have
$C_{S^1,\; \ast}(\Psi_{a, K}, \;S^1\cdot u_K)\cong C_{S^1,\; q}(F_{a, K}, \;S^1\cdot x)$. By Proposition 3.2,
the module $C_{S^1,\; q}(F_{a, K}, \;S^1\cdot x)$ is independent of the choice of the Hamiltonian function $H_a$
whenever $H_a$ satisfies conditions in Proposition 2.5. Hence in order to compute the critical modules, we can
choose $\Psi_{a, K}$ with $H_a$ being positively homogeneous of degree $\aa=\aa_a$ near the image set of every
nonzero solution $x$ of (\ref{2.2}) which corresponds to some closed characteristic $(\tau,y)$ with period $\tau$
being strictly less than $aT$.

In other words, for a given $a>0$, we choose $\vth\in (0,1)$ first such that
$[aT\vth, aT(1-\vth)]\supset \per(\Sigma)\cap (0,aT)$ holds by the definition of the set $\per(\Sg)$ and the assumption
(F). Then we choose $\aa=\aa_a\in (1,¡¢£¬2)$ sufficiently close to $2$ by (iii) of Lemma 2.2 such that $\vf_a(t)=ct^\aa$
for some constant $c>0$ and $\aa\in(1,\,2)$ whenever $\frac{\vf_a^\prime(t)}{t}\in [\vth, 1-\vth]$. Now we suppose that
$\vf_a$ satisfies (iii) of Lemma 2.2.

Now we consider iterations of critical points of $\Psi_{a, K}$. Suppose $u_K\neq 0$ is a critical point of
$\Psi_{a,K}$ with $\mul(u_K)=m$, and $g_K$ is the critical point of $\psi_{a,K}$ corresponding to $u_K$. By
Proposition 2.5 and Lemma 2.7, we have $u_K=-J\dot x+Kx$ with $x$ being a solution of (\ref{2.2}) and
$x=\rho y(\frac{\tau t}{T})$ with $\frac{\vf_a^\prime(\rho)}{\rho}=\frac{\tau}{aT}$. Moreover, $(\tau, y)$ is a closed
characteristic on $\Sigma$ with minimal period $\frac{\tau}{m}$. For any $p\in\N$ satisfying $p\tau<aT$, we choose $K$
such that $pK\notin \frac{2\pi}{T}\Z$, then the $p$th iteration $u_{pK}^p$ of $u_K$ is given by $-J\dot x^p+pKx^p$,
where $x^p$ is the unique solution of (2.2) corresponding to $(p\tau, y)$ and is a critical point of $F_{a, pK}$, that
is, $u_{pK}^p$ is the critical point of $\Psi_{a, pK}$ corresponding to $x^p$. Hence we have
\bea
&& x(t)=\left(\frac{\tau}{c\alpha a}\right)^\frac{1}{\alpha-2}y(\tau t), \quad
      x^p(t)=\left(\frac{p\tau}{c\alpha a}\right)^\frac{1}{\alpha-2}y(p\tau t)=p^{\frac{1}{\alpha-2}}x(pt), \nn\\
&& u_K(t)=-J\dot {x}(t)+Kx(t), \quad u_{pK}^p(t)=-J\dot x^p(t)+pKx^p(t)=p^{\frac{\alpha-1}{\alpha-2}}u_K(pt). \nn\eea
We define the $p$th iteration $\phi^p$ on $L^{2}(\R/(T\Z); {\bf R}^{2n})$ by
\be \phi^p: v_K(t)\mapsto v^p_{pK}(t)\equiv p^\frac{\alpha-1}{\alpha-2}v_K(pt). \lb{4.3}\ee
Then there exist a $w\in W^{1, 2}(\R/(T\Z); {\bf R}^{2n})$ such that
\be  v_K(t)=-J\dot{w}(t)+Kw(t),\quad v_{pK}^p(t)=-J\dot{w}^p(t)+pKw^p(t),
              \quad w^p(t)=p^{\frac{1}{\alpha-2}}w(pt).  \lb{4.4}\ee
By definition, we have
\bea
\Psi_{a,K}(v_K) &=& \int_0^T{\left[\frac{1}{2}(J\dot{w}-K w,w)+H_{a,K}^*(-J\dot{w}+K w)\right]dt}  \lb{4.5}\\
\Psi_{a, pK}(v^p_{pK}) &=& \int_0^T{\left[\frac{1}{2}(J\dot{w}^p-pKw^p, w^p)+H_{a,pK}^*(-J\dot{w}^p+pKw^p)\right]dt}\nn\\
&=&\frac{1}{2}p^{\frac{\alpha}{\alpha-2}}\int_0^T{(J\dot{w}(pt)-K w(pt),w(pt))dt}\nn\\
&&\quad + \int_0^T{H_{a,pK}^*(p^{\frac{\alpha-1}{\alpha-2}}(-J\dot{w}(pt)+Kw(pt)))dt}.   \lb{4.6}\eea
Let $\xi(t)=H_{a, K}^{*\prime}(-J\dot{w}(t)+Kw(t))$, then $H_{a, K}^\prime(\xi(t))=-J\dot{w}(t)+Kw(t)$. Note that when
$v_K$ belongs to a small $L^\infty$-neighborhood of $u_K$, $\xi$ belongs to a small $L^\infty$-neighborhood of $x$. In
the following, we suppose that $v_K$ belongs to a small $L^\infty$-neighborhood of $u_K$.

Since $H_a$ is positively homogeneous of degree $\aa$ near the image set of $x$, we have
$H_a^\prime(p^{\frac{1}{\aa-2}}\xi(pt))=p^{\frac{\aa-1}{\aa-2}}H_a^\prime(\xi(pt))$,
$H_{a}(p^{\frac{1}{\aa-2}}\xi(pt))=p^{\frac{\aa}{\aa-2}}H_a(\xi(pt))$. Thus, we get
\bea H_{a,pK}^\prime(p^{\frac{1}{\aa-2}}\xi(pt))
&=& H_a^\prime(p^{\frac{1}{\aa-2}}\xi(pt))+pKp^{\frac{1}{\aa-2}}\xi(pt) \nn\\
&=& p^{\frac{\aa-1}{\aa-2}}H_{a,K}^\prime(\xi(pt))  \nn\\
&=& p^{\frac{\alpha-1}{\alpha-2}}(-J\dot{w}(pt)+Kw(pt)). \nn\eea
So it follows that
\bea
H_{a,pK}^*(p^{\frac{\alpha-1}{\alpha-2}}(-J\dot{w}(pt)+Kw(pt)))
&=& H_{a, pK}^\prime(p^{\frac{1}{\aa-2}}\xi(pt))\cdot p^{\frac{1}{\aa-2}}\xi(pt)-H_{a,pK}(p^{\frac{1}{\aa-2}}\xi(pt)) \nn\\
&=& p^{\frac{\aa}{\aa-2}}H_{a,K}^\prime(\xi(pt))\cdot \xi(pt)-p^{\frac{\aa}{\aa-2}}H_a(\xi(pt))-pK|p^{\frac{1}{\aa-2}}\xi(pt)|^2 \nn\\
&=& p^{\frac{\aa}{\aa-2}}(H_{a,K}^\prime(\xi(pt))\cdot\xi(pt)-H_{a,K}(\xi(pt))) \nn\\
&=& p^{\frac{\aa}{\aa-2}}H_{a,K}^*(H_{a,K}^\prime(\xi(pt))) \nn\\
&=& p^{\frac{\aa}{\aa-2}}H_{a,K}^*(-J\dot{w}(pt)+Kw(pt))   \lb{4.7}\eea
Combining (\ref{4.5})-(\ref{4.7}), we obtain
\bea  \Psi_{a,pK}(v^p_{pK})
&=& p^{\frac{\alpha}{\alpha-2}}\int_0^T{\frac{1}{2}(J\dot{w}(pt)-Kw(pt), w(pt))+H_{a,K}^*(-J\dot{w}(pt)+Kw(pt))dt} \nn\\
&=& p^{\frac{\alpha}{\alpha-2}}\int_0^{pT}{\frac{1}{2}(J\dot{w}(t)-Kw(t), w(t))+H_{a,K}^*(-J\dot{w}(t)+Kw(t))d(\frac{t}{p})} \nn\\
&=& p^{\frac{\alpha}{\alpha-2}}\Psi_{a,K}(v_K).  \lb{4.8}\eea
By direct computation, we obtain
$$ \Psi_{a,pK}^\prime(v_{pK}^p) = -w^p+H^{*\prime}_{a,pK}(v_{pK}^p)=-p^{\frac{1}{\alpha-2}}w(pt)+p^{\frac{1}{\aa-2}}\xi(pt), $$
specially,
\be \Psi_{a,K}^\prime(v_{K}) = -w+H^{*\prime}_{a,K}(v_{K})=-w(t)+\xi(t). \lb{4.9}\ee
Hence, we have
\be \Psi_{a, pK}^\prime(\phi^p(v_K)) = \Psi_{a,pK}^\prime(v_{pK}^p) = p^{-1}\phi^p(\Psi_{a, K}^\prime(v_{K})). \lb{4.10}\ee

Applying it to Lemma 2.10, for $g\in DN(S^1\cdot g_K)$, noticing that $g$ and $h_K(g)\in C^0(\R/\Z; \R^{2n})$,
and $h_K: G_K\rightarrow G^\bot_K$ is continuous in the $C^0$-topology (cf. Page 628 of \cite{Vit1}), where we write
$h_K$ and $G_K$ for $h_a$ and $G$ respectively to indicate their dependence on $K$, we have that $g+h_K(g)$
belongs to a small $L^\infty$-neighborhood of $u_K$ when the radius $\varrho>0$ of the ball $DN(S^1\cdot g_K)$ is small
enough and
\be \frac{\partial}{\pt h}\Psi_{a,pK}(\phi^p(g)+h_{pK}(\phi^p(g)))=0=\frac{\pt}{\pt h}\Psi_{a, pK}(\phi^p(g)+\phi^p(h_{K}(g))),
          \lb{4.11}\ee
where we choose $G_{pK}$ in Lemma 2.10 for $\Psi_{a, pK}$ such that $G_{pK}\supseteq \phi^p(G_K)$, in fact, in (\ref{2.18}) we can
choose the same $\omega>0$ for both $K$ and $pK$, and let $G_{pK}$ be the subspace of $L^2(\R/\Z, \R^{2n})$ generated
by the eigenvectors of $-M_{pK}$ whose eigenvalues are less than $-\frac{\omega}{2p}$; that is,
$$ G_{pK} = \span\{e^{-JLt}x_0\mid -\frac{1}{L+pK}<-\frac{\omega}{2p}, L\in\frac{2\pi}{T}\Z, x_0\in\R^{2n}\},  $$
and then
$$ \phi^p(G_K) = \span\{e^{-JpLt}x_0\mid -\frac{1}{L+K}<-\frac{\omega}{2}, L\in\frac{2\pi}{T}\Z, x_0\in\R^{2n}\}
                 \subseteq G_{pK}. $$
Hence $h_{pK}(\phi^p(g))=\phi^p(h_K(g))$ holds. This together with (\ref{4.8}) yields
\bea \psi_{a, pK}(\phi^p(g))
&=& \Psi_{a,pK}(\phi^p(g)+h_{pK}(\phi^p(g))) =\Psi_{a,pK}(\phi^p(g)+\phi^p(h_K(g))) \nn\\
&=&p^{\frac{\alpha}{\alpha-2}}\Psi_{a,K}(g+h_K(g))=p^{\frac{\alpha}{\alpha-2}}\psi_{a, K}(g).  \lb{4.12}\eea

We define a new inner product $\<\cdot,\cdot\>_p$ on $L^2(\R/\Z, \R^{2n})$ by
\be \<v,w\>_p = p^\frac{2(\alpha-1)}{2-\alpha}\<v, w\>.   \lb{4.13}\ee
Then $\phi^p: DN(g_K)\rightarrow DN(g_{pK}^p)$ is an isometry from the standard inner product to the above one,
where $g_{pK}^p=\phi^p(g_K)$ is the critical point of $\psi_{a,pK}$ corresponding to $u_{pK}^p$ and the radii of
the two normal disk bundles are suitably chosen. Clearly $\phi^p(DN(g_K))$ consists of points in $DN(g_{pK}^p)$
which are fixed by the $\Z_p$-action. Since the $\Z_p$-action on $DN(g_{pK}^p)$ are isometries and
$f\equiv\psi_{a, pK}|_{DN(g_{pK}^p)}$ is $\Z_p$-invariant, we have
\be f^{\prime\prime}(g) =\left(\matrix{(f|_{\phi^p(DN(g_K))})^{\prime\prime} \quad 0\cr \qquad
             0\qquad\qquad\;\;\ast}\right),\quad \forall g\in \phi^p(DN(g_K)).  \lb{4.14}\ee
Moreover, we have
\be f^\prime(g)=(f|_{\phi^p(DN(g_K))})^{\prime},\quad \forall g\in \phi^p(D N(g_K)).  \lb{4.15}\ee
Now we can apply the results by D. Gromoll and W. Meyer \cite{GrM1} to the manifold $DN(g_{pK}^p)$ with $g_{pK}^p$
as its unique critical point. Then $\mul(g_{pK}^p)=pm$ is the multiplicity of $g_{pK}^p$ and the isotropy group
$\Z_{pm}\subseteq S^1$ of $g_{pK}^p$ acts on $DN(g_{pK}^p)$ by isometries. According to Lemma 1 of \cite{GrM1}, we
have a $\Z_{pm}$-invariant decomposition of $T_{g_{pK}^p}(DN(g_{pK}^p))$
\be  T_{g_{pK}^p}(DN(g_{pK}^p)) = V^+\oplus V^-\oplus V^0 = \{(x_+, x_-, x_0)\}  \lb{4.16}\ee
with $\dim V^-=i(g_{pK}^p)=i_{pK}(u_{pK}^p)$, $\dim V^0=\nu(g_{pK}^p)-1=\nu_{pK}(u_{pK}^p)-1$ (cf. Lemma 2.10(iii)),
and a $\Z_{pm}$-invariant neighborhood $B=B_+\times B_-\times B_0$ for $0$ in $T_{g_{pK}^p}(DN(g_{pK}^p))$ together
with two $\Z_{pm}$-invariant diffeomorphisms
$$ \Phi : B=B_+\times B_-\times B_0\to \Phi(B_+\times B_-\times B_0)\subset DN(g_{pK}^p),  $$
and
$$ \eta : B_0\to W(g_{pK}^p)\equiv\eta(B_0)\subset DN(g_{pK}^p), $$
and $\Phi(0)=\eta(0)=g_{pK}^p$, such that
\be \psi_{a,pK}\circ\Phi(x_+,x_-,x_0)=|x_+|^2 - |x_-|^2 + \psi_{a,pK}\circ\eta(x_0), \lb{4.17}\ee
with $d(\psi_{a, pK}\circ \eta)(0)=d^2(\psi_{a, pK}\circ\eta)(0)=0$. As usual, we call $W(g_{pK}^p)$ a local
characteristic manifold, and $U(g_{pK}^p)=B_-$ a local negative disk at $g_{pK}^p$. By the proof of Lemma 1 of
\cite{GrM1}, $W(g_{pK}^p)$ and $U(g_{pK}^p)$ are $\Z_{pm}$-invariant. It follows from (\ref{4.17}) that $g_{pK}^p$
is an isolated critical point of $\psi_{a, pK}|_{DN(g_{pK}^p)}$. Then as in Lemma 6.4 of \cite{Rad2}, we have
\bea
&& H_\ast(\wtd{\Lambda}_{a, pK}(g_{pK}^p)\cap DN(g_{pK}^p),\;
     (\wtd{\Lambda}_{a, pK}(g_{pK}^p)\setminus\{g_{pK}^p\})\cap DN(g_{pK}^p))\nn\\
&&\quad = \bigoplus_{q\in\Z}H_q (U(g_{pK}^p),U(g_{pK}^p)\setminus\{g_{pK}^p\}) \nn\\
&&\qquad\qquad\qquad \otimes H_{\ast-q}(W(g_{pK}^p)\cap \wtd{\Lambda}_{a,pK}(g_{pK}^p),
       (W(g_{pK}^p)\setminus\{g_{pK}^p\})\cap \wtd{\Lambda}_{a,pK}(g_{pK}^p)), \lb{4.18}\eea
where
\be  H_q(U(g_{pK}^p),U(g_{pK}^p)\setminus\{g_{pK}^p\} )
    = \left\{\matrix{\Q, & {\rm if\;}q=i_{pK}(u_{pK}^p),  \cr
                      0, & {\rm otherwise}. \cr}\right.    \lb{4.19}\ee

Now we have the following proposition.

{\bf Proposition 4.2.} {\it For any $p\in\N$, we choose $K$ such that $pK\notin \frac{2\pi}{T}\Z$. Let $u_K\neq 0$
be a critical point of $\Psi_{a, K}$ with $\mul(u_K)=1$, $u_K=-J\dot x+Kx$ with $x$ being a critical point of
$F_{a, K}$. Then for all $q\in\Z$, we have
\bea
&& C_{S^1,\; q}(\Psi_{a,pK},\;S^1\cdot u_{pK}^p) \nn\\
&&\quad\cong \left(\frac{}{}H_{q-i_{pK}(u_{pK}^p)}(W(g_{pK}^p)\cap \wtd{\Lambda}_{a,pK}(g_{pK}^p),(W(g_{pK}^p)
                 \setminus\{g_{pK}^p\})\cap \wtd{\Lambda}_{a,pK}(g_{pK}^p))\right)^{\beta(x^p)\Z_p}, \qquad\quad \lb{4.20}\eea
where $\beta(x^p)=(-1)^{i_{pK}(u_{pK}^p)-i_K(u_K)}=(-1)^{i^v(x^p)-i^v(x)}$. In particular, if $u_{pK}^p$ is
non-degenerate, i.e., $\nu_{pK}(u_{pK}^p)=1$, then}
\be C_{S^1,\; q}(\Psi_{a,pK},\;S^1\cdot u_{pK}^p)
    = \left\{\matrix{\Q, & {\rm if\;}q=i_{pK}(u_{pK}^p)\;{\rm and\;}\beta(x^p)=1,  \cr
                      0, & {\rm otherwise}. \cr}\right. \lb{4.21}\ee

{\bf Proof.} Suppose $\theta$ is a generator of the linearized $\Z_p$-action on $U(g_{pK}^p)$. Then $\theta(\xi)=\xi$
if and only if $\xi\in T_{g_{pK}^p}(\phi^p(DN(g_K)))$. Hence it follows from (\ref{4.12}) and (\ref{4.14}) that
$\xi=(\phi^p)_\ast(\xi^\prime)$ for a unique $\xi^\prime\in T_{g_K}(DN(g_K))^-$. Hence the proof of Satz 6.11 in
\cite{Rad2}, Proposition 2.8 in \cite{BaL1} or Proposition 3.10 in \cite{WHL1} yield this proposition. Note that
$i_{pK}(u_{pK}^p)=2n([pKT/{2\pi}]+1)+i^v(x^p)$ and $i_{K}(u_{K})=2n([KT/{2\pi}]+1)+i^v(x)$ follow from
(\ref{2.14}) and (\ref{2.15}). \hfill\hb

{\bf Definition 4.3.} {\it For any $p\in\N$, we choose $K$ such that $pK\notin \frac{2\pi}{T}\Z$. Let $u_K\neq 0$ be
a critical point of $\Psi_{a,K}$ with $\mul(u_K)=1$, $u_K=-J\dot x+Kx$ with $x$ being a critical point of $F_{a, K}$.
Then for all $l\in\Z$, let
\bea
k_{l,\pm 1}(u_{pK}^p) &=& \dim\left(\frac{}{}H_{l}(W(g_{pK}^p)\cap
  \wtd{\Lambda}_{a,pK}(g_{pK}^p),(W(g_{pK}^p)\bs\{g_{pK}^p\})\cap \wtd{\Lambda}_{a,pK}(g_{pK}^p))\right)^{\pm\Z_p},
           \quad  \lb{4.22}\\
k_l(u_{pK}^p) &=& \dim\left(\frac{}{}H_{l}(W(g_{pK}^p)\cap
  \wtd{\Lambda}_{a,pK}(g_{pK}^p),(W(g_{pK}^p)\bs\{g_{pK}^p\})\cap \wtd{\Lambda}_{a,pK}(g_{pK}^p))\right)^{\beta(x^p)\Z_p}.
           \qquad\quad   \lb{4.23}\eea
Here $k_l(u_{pK}^p)$'s are called critical type numbers of $u_{pK}^p$. }

{\bf Remark 4.4.} (i) Since
\bea C_{S^1,\;l+i_{pK}(u_{pK}^p)}(\Psi_{a,pK}, \;S^1\cdot u_{pK}^p)
&\cong& C_{S^1,\; l+i_{pK}(x^p)}(F_{a,pK}, \;S^1\cdot x^p)  \nn\\
&\cong& C_{S^1,\; l+d(pK)+i^v(x^p)}(F_{a,pK}, \;S^1\cdot x^p),  \nn\eea
by Theorem 3.3, we obtain that $k_l(u_{pK}^p)$ is independent of the choice of $K$ and denote it by $k_l(x^p)$,
here $k_l(x^p)$'s are called critical type numbers of $x^p$.

(ii) By Proposition 2.11, we have $k_{l, \pm 1}(u_{pK}^p)=0$ if $l\notin [0, 2n-2]$.

Similar to Section 7.1 of \cite{Rad2}, Theorem 2.11 of \cite{BaL1}, or Lemma 3.12 of \cite{WHL1}, we have

{\bf Lemma 4.5.} {\it Let $u_K\neq 0$ be a critical point of $\Psi_{a, K}$ with $\mul(u_K)=1$. Suppose
$\nu_{mK}(u_{mK}^m)=\nu_{pmK}(u_{pmK}^{pm})$ for some $m, p\in\N$. Then we have
$k_{l, \pm 1}(u_{mK}^m)=k_{l, \pm 1}(u_{pmK}^{pm})$ for all $l\in\Z$.}

{\bf Proof.} Let $\phi^p: DN(g_{mK}^m)\rightarrow DN(g_{pmK}^{pm})$ be the $p$th iteration map. By (\ref{4.13}),
$\phi^p$ is an isometry under the modified metric. Hence by (\ref{4.12}),  we have
\be  \nu_{mK}(u_{mK}^m)-1 = \dim\ker((\psi_{a,mK}|_{DN(g_{mK}^m)})^{\prime\prime}-I)
        = \dim\ker((\psi_{a,pmK}|_{\phi^p(DN(g_{mK}^m))})^{\prime\prime}-I).  \lb{4.24}\ee
Thus by (\ref{4.14}) and the assumption $\nu_{mK}(u_{mK}^m)=\nu_{pmK}(u_{pmK}^{pm})$, we have that
$T_{g_{pmK}^{pm}}(\phi^p(DN(g_{mK}^m)))$ contains the null space of the Hessian of $\psi_{a, pmK}|_{DN(g_{pmK}^{pm})}$.
Now by (\ref{4.15}), we can use Lemma 7 of \cite{GrM1} to obtain that $\phi^p(W(g_{mK}^m))\equiv W(g_{pmK}^{pm})$ is
a characteristic manifold of $\psi_{a,pmK}|_{DN(g_{pmK}^{pm})}$, where $W(g_{mK}^m)$ is a characteristic manifold of
$\psi_{a,mK}|_{DN(g_{mK}^m)}$. By (\ref{4.12}), we have
\bea \phi^p:
&&(W(g_{mK}^m)\cap\wtd{\Lambda}_{a,mK}(g_{mK}^m),(W(g_{mK}^m)\setminus\{g_{mK}^m\})\cap\wtd{\Lambda}_{a,pK}(g_{mK}^m)) \nn\\
&&\qquad \rightarrow (W(g_{pmK}^{pm})\cap\wtd{\Lambda}_{a,pmK}(g_{pmK}^{pm}),
       (W(g_{pmK}^{pm})\setminus\{g_{pmK}^{pm}\})\cap\wtd{\Lambda}_{a, pmK}(g_{pmK}^{pm}))  \nn\eea
is a homeomorphism. Suppose $\theta$ and $\theta_p$ generate the $\Z_m$ and $\Z_{pm}$ action on $W(g_{mK}^m)$ and
$W(g_{pmK}^{pm})$ respectively. Then clearly $\phi^p\circ \theta=\theta_p\circ \phi^p$ holds and it implies
\bea
&& H_\ast(W(g_{mK}^m)\cap\wtd{\Lambda}_{a,mK}(g_{mK}^m),(W(g_{mK}^m)\setminus\{g_{mK}^m\})
                 \cap \wtd{\Lambda}_{a,pK}(g_{mK}^m))^{\pm\Z_m} \nn\\
&&\qquad \cong H_*(W(g_{pmK}^{pm})\cap\wtd{\Lambda}_{a,pmK}(g_{pmK}^{pm}),(W(g_{pmK}^{pm})\setminus\{g_{pmK}^{pm}\})
                  \cap\wtd{\Lambda}_{a, pmK}(g_{pmK}^{pm}))^{\pm\Z_{pm}}. \nn\eea
Therefore our lemma follows. \hfill\hb

{\bf Proposition 4.6.} {\it Let $x\neq 0$ be a critical point of $F_{a,K}$ with $\mul(x)=1$ corresponding to a
critical point $u_K$ of $\Psi_{a, K}$. Then there exists a minimal $K(x)\in \N$ such that
\bea
&& \nu^v(x^{p+K(x)})=\nu^v(x^p),\quad i^v(x^{p+K(x)})-i^v(x^p)\in 2\Z,  \qquad\forall p\in \N,  \lb{4.25}\\
&& k_l(x^{p+K(x)})=k_l(x^p), \qquad\forall p\in \N,\;l\in\Z. \lb{4.26}\eea
We call $K(x)$ the minimal period of critical modules of iterations of the functional $F_{a, K}$ at $x$. }

{\bf Proof.} As in the proof of Proposition 2.11, we denote by $R(t)$ the fundamental solution of (\ref{2.16}).
Then by Section 2 and Theorem 2.1 of \cite{HuL1}, we have $i^v(x^p)=i(x,p)-n$ and $\nu^v(x^p)=\nu(x,p)$ for all
$p\in\N$, where $(i(x,p), \nu(x,p))$ are index and nullity defined by C. Conley, E. Zehnder and Y. Long
(cf. \cite{CoZ1}, \cite{LZe1}, \cite{Lon1}-\cite{Lon4}). Hence we have $\nu^v(x^p)=\dim\ker (R(1)^p-I_{2n})$.
Denote by $\lm_i=\exp(\pm 2\pi\frac{r_i}{s_i})$ the eigenvalues of $R(1)$ possessing rotation angles which are
rational multiples of $\pi$ with $r_i$ and $s_i\in\N$ and $(r_i,s_i)=1$ for $1\le i\le q$. Let $K(x)$ be twice
of the least common multiple of $s_1,\ldots, s_q$. Then (\ref{4.25}) holds. Note that the later conclusion in
(\ref{4.25}) follows from Theorem 9.3.4 of \cite{Lon4}.

In order to prove (\ref{4.26}), it suffices to show
\be  k_l(x^{m+qK(x)})=k_l(x^m), \qquad \forall q\in \N, \; l\in\Z, \; 1\le m\le K(x).  \lb{4.27}\ee

In fact, assume that (\ref{4.27}) is proved. Note that (\ref{4.26}) follows from (\ref{4.27}) with $q=1$
directly when $p\le K(x)$. When $p>K(x)$, we write $p=m+qK(x)$ for some $q\in\N$ and $1\le m\le K(x)$. Then by
(\ref{4.27}) we obtain
$$ k_l(x^{p+K(x)})=k_l(x^{m+(q+1)K(x)})=k_l(x^m)=k_l(x^{m+qK(x)})=k_l(x^p), $$
i.e., (\ref{4.26}) holds.

To prove (\ref{4.27}), we fix an integer $m\in [1,K(x)]$. Let
$$ A=\{s_i\in\{s_1,\ldots, s_q\}\;|\;\;s_i\;{\rm is\;a\;factor\;of\;}m\},  $$
and let $m_1$ be the least common multiple of elements in $A$. Hence we have $m=m_1m_2$ for some $m_2\in \N$
and $\nu_{mK}(u_{mK}^{m})=\nu^v(x^m)=\nu^v(x^{m_1})=\nu_{m_1K}(u_{m_1K}^{m_1})$. Thus by Remark 4.4 (i) and
Lemma 4.5, we have $k_l(x^m)=k_{l,\beta(x^m)}(u_{mK}^{m})=k_{l,\beta(x^m)}(u_{m_1K}^{m_1})$. Since
$m+pK(u)=m_1m_3$ for some $m_3\in\N$, we have by Remark 4.4 (i) and Lemma 4.5 that
$k_l(x^{m+pK(x)})=k_{l, \beta(x^{m+pK(x)})}(u_{m_1K}^{m_1})$. By (\ref{4.25}), we obtain
$\beta(x^{m+pK(x)})=\beta(x^m)$, and then (\ref{4.27}) is proved. This completes the proof. \hfill\hb

Note that the above Proposition 4.6 could be established also without forcing the Hamiltonian
to be homogeneous near its critical points. In fact, by Proposition 3.2, it holds for any Hamiltonian
defined by Proposition 2.5.

In the following, Let $F_{a,K}$ be the functional defined by (\ref{2.6}) with $H_a$ satisfying Proposition 2.5,
we do not require $\wtd{H}_a$ to be homogeneous anymore.

{\bf Definition 4.7.} {\it Suppose the condition (F) at the beginning of Section 2 holds. For every closed
characteristic $(\tau,y)$ on $\Sigma$, let $aT>\tau$ and choose $\vf_a$ to satisfy (i)-(ii) of Lemma 2.2.
Determine $\rho$ uniquely by $\frac{\vf_a'(\rho)}{\rho}=\frac{\tau}{aT}$. Let $x=\rho y(\frac{\tau t}{T})$.
Then we define the index $i(\tau,y)$ and nullity $\nu(\tau,y)$ of $(\tau,y)$ by
$$ i(\tau,y)=i^v(x), \qquad \nu(\tau,y)=\nu^v(x). $$
Then the mean index of $(\tau, y)$ is defined by }
\be \hat i(\tau,y) = \lim_{m\rightarrow\infty}\frac{i(m\tau, y)}{m}.  \lb{4.28}\ee

Note that by Proposition 2.11, the index and nullity are well defined and is independent of the
choice of $aT>\tau$ and $\vf_a$ satisfying (i)-(ii) of Lemma 2.2.

For a prime closed characteristic $(\tau, y)$ on $\Sigma$, we denote simply by $y^m\equiv(m\tau, y)$ for $m\in\N$.
By Proposition 3.2, we can define the critical type numbers $k_l(y^m)$ of $y^m$ to be $k_l(x^m)$, where $x^m$
is the critical point of $F_{a, K}$ corresponding to $y^m$. We also define $K(y)=K(x)$, where $K(x)\in \N$
is given by Proposition 4.6. Suppose $\Nn$ is an $S^1$-invariant open neighborhood of $S^1\cdot x^m$ such that
$\crit(F_{a, K})\cap(X_{a,K}(x^m)\cap\Nn)=S^1\cdot x^m$. Then we make the following definition

{\bf Definition 4.8.} {\it The Euler characteristic $\chi(y^m)$ of $y^m$ is defined by
\bea \chi(y^m)
&\equiv& \chi((X_{a,K}(x^m)\cap\Nn)_{S^1},\; ((X_{a,K}(x^m)\setminus S^1\cdot x^m)\cap\Nn)_{S^1}) \nn\\
&\equiv& \sum_{q=0}^{\infty}(-1)^q\dim C_{S^1,\;q}(F_{a,K},\;S^1\cdot x^m). \lb{4.29}\eea
Here $\chi(A,B)$ denotes the usual Euler characteristic of the space pair $(A,B)$.
The average Euler characteristic $\hat\chi(y)$ of $y$ is defined by}
\be \hat{\chi}(y)=\lim_{N\to\infty}\frac{1}{N}\sum_{1\le m\le N}\chi(y^m). \lb{4.30}\ee

Note that by Proposition 3.2 and Theorem 3.3, $ \chi(y^m)$ is well defined and is independent of the
choice of $a$ and $K$. In fact, by Remark 4.4 (i), we have
\be  \chi(y^m)=\sum_{l=0}^{2n-2}(-1)^{i(y^m)+l}k_l(y^m).  \lb{4.31}\ee
The following remark shows that $\hat\chi(y)$ is well-defined and is a rational number.

{\bf Remark 4.9.} By (\ref{4.25}), (\ref{4.31}) and Proposition 4.6, we have
\bea \hat\chi(y)
&=&\lim_{N\rightarrow\infty}\frac{1}{N}  \sum_{1\le m\le N\atop 0\le l\le 2n-2}(-1)^{i(y^m)+l}k_l(y^m) \nn\\
&=&\lim_{s\rightarrow\infty}\frac{1}{sK(y)}\sum_{1\le m\le K(y),\; 0\le l\le 2n-2\atop 0\le p< s}
              (-1)^{i(y^{pK(y)+m})+l}k_l(y^{pK(y)+m}) \nn\\
&=&\frac{1}{K(y)}\sum_{1\le m\le K(y)\atop 0\le l\le 2n-2}(-1)^{i(y^{m})+l}k_l(y^{m}).  \lb{4.32}\eea
Therefore $\hat\chi(y)$ is well defined and is a rational number. In particular, if all $y^m$s are
non-degenerate, then $\nu(y^m)=1$ for all $m\in\N$. Hence the proof of Proposition 4.6 yields $K(y)=2$.
By (\ref{4.21}), we have
$$ k_l(y^m)
    = \left\{\matrix{1, & {\rm if\;\;} i(y^m)-i(y)\in 2\Z \quad {\rm and} \quad l=0  \cr
                     0, & {\rm otherwise}. \cr}\right.  $$
Hence (\ref{4.32}) implies
\be \hat\chi(y)
    = \left\{\matrix{(-1)^{i(y)}, & {\rm if\;\;} i(y^2)-i(y)\in 2\Z,  \cr
           \frac{(-1)^{i(y)}}{2}, & {\rm otherwise}. \cr}\right.  \lb{4.33}\ee

{\bf Remark 4.10.} Note that $k_l(y^m)=0$ for $l\notin [0, \nu(y^m)-1]$ and it can take only values $0$
or $1$ when $l=0$ or $l=\nu(y^m)-1$. Moreover, the following facts are useful (cf. Lemma 3.11 of
\cite{BaL1}, Remark 3.17 of \cite{WHL1}, \cite{Cha1} and \cite{MaW1}):

(i) $k_0(y^m)=1$ implies $k_l(y^m)=0$ for $1\le l\le \nu(y^m)-1$.

(ii) $k_{\nu(y^m)-1}(y^m)=1$ implies $k_l(y^m)=0$ for $0\le l\le \nu(y^m)-2$.

(iii) $k_l(y^m)\ge 1$ for some $1\le l\le \nu(y^m)-2$ implies $k_0(y^m)=k_{\nu(y^m)-1}(y^m)=0$.

(iv) In particular, only one of the $k_l(y^m)$s for $0\le l\le \nu(y^m)-1$ can be non-zero when $\nu(y^m)\le 3$.

\setcounter{equation}{0}
\section{Contribution of the origin}

In section 3 and 4, we studied nonzero critical points of $F_{a,K}$, now we need to study the contribution
of the origin to the Morse series of the functional $F_{a,K}$ on $W^{1,2}(\R/\Z;\R^{2n})$. Theorem 7.1 of
\cite{Vit1} was given under the condition that all the closed characteristics together with their iterations
are non-degenerate, however, by a modification of the proof, we obtain a degenerate version in the following.

{\bf Theorem 5.1.} {\it Fix an $a>0$ such that $\per(\Sg)\cap (0,aT)\not=\emptyset$. Then there exists an
$\vep_0>0$ small enough such that for any $\vep\in (0,\vep_0]$ we have
\be H_{S^1,\;q+d(K)}(X_{a, K}^{\vep}, \;X_{a, K}^{-\vep}) = 0,  \quad \forall q\in \mathring {I}, \lb{5.1}\ee
if $I$ is an interval of $\Z$ such that $I\cap [i(\tau, y), i(\tau, y)+\nu(\tau, y)-1]=\emptyset$ for all
closed characteristics $(\tau,\, y)$ on $\Sigma$ with $\tau\ge aT$.}

{\bf Proof.} By Proposition 2.8, there exists an $\vep_0>0$ such that there are no critical values of
$F_{a,K}$ in the interval $[-\vep_0,\; \vep_0]$ except $0$. Hence we have
$$ H_{S^1,\; q+d(K)}(X_{a, K}^{\vep}, \;X_{a, K}^{-\vep}) \cong
  H_{S^1,\;q+d(K)}(X_{a, K}^{\vep_0}, \;X_{a, K}^{-\vep_0}), \quad \forall q\in\Z,\; \vep\in(0,\,\vep_0]. $$
In the following we assume $\vep\in(0,\,\vep_0]$.

Note that by the same proof of Proposition 3.2, $H_{S^1,\; q+d(K)}(X_{a, K}^{\vep}, \;X_{a, K}^{-\vep})$ is
independent of the choice of $\vf_a$ in $\wtd{H}_a(x)=a\vf_a(j(x))$ which satisfies (i) of Lemma 2.2. Hence
we can choose $a\vf_a\equiv\phi$, where $\phi$ is defined as in Lemma 2.2 of \cite{Vit1}. Since the homology
in (\ref{5.1}) depends only on the restriction of $H_a$ to a neighborhood of the origin, as in the beginning
of Section 7 of \cite{Vit1}, we assume $H_a$ to be homogeneous of degree two everywhere. Now we can make
some modifications of the proof of Theorem 7.1 of \cite{Vit1} to complete our proof.

There are only four palaces in the proof of Theorem 7.1 of \cite{Vit1} where the non-degenerate
condition is used, i.e. the use of Proposition 1 of Appendix 1 (In fact, it should be Proposition 3
of Appendix 1) in Page 640 and Page 642, the use of Proposition 2 (i) of Appendix 1 in Page 644,
and the use of Proposition 3 of Appendix 1 in Page 648.

However, Proposition 2 (i), Proposition 3 (a) of Appendix 1 of \cite{Vit1} work for the degenerate
case, so we only need to modify the proof in Page 648 of \cite{Vit1}. By assumption,
$I\cap [i(\tau, y), i(\tau, y)+\nu(\tau, y)-1]=\emptyset$ for all closed characteristics $(\tau,\, y)$ on
$\Sigma$ with $\tau\ge aT$, then by a degenerate version of Proposition 3 (b) of Appendix 1 of \cite{Vit1},
we complete our modification. More precisely, in the proof of Proposition 3 (b) of Appendix 1 of \cite{Vit1}
(we only consider Case (i)), we let $\hat{\tau}$ be the restriction of the map $\tau$ on $\partial X$, then
\be C_*(\hat{\tau}, (x_0,t_0))\cong C_*(f_{t_0}, x_0).  \lb{5.2}\ee
In fact, since $\partial f_t(x_0)/{\partial t}<0$, by the implicit functional theorem, we have that in a
small enough neighborhood $U$ of $x_0$, there is an unique continuous function $t_x$ such that
$(x, t_x)\in \partial X$ and $t_{x_0}=t_0$ for $x\in U$. By the fact that $\partial f_t(x_0)/{\partial t}<0$,
we obtain that $t_x\leq t_0$ for $x\in U$ if and only if $f_{t_0}(x)\leq 0$. Thus
$\{(x, t_x)\mid\hat{\tau}(x, t_x)\leq t_0, x\in U\}$ is homotopy equivalent with
$\{x\in U\mid f_{t_0}(x)\leq 0\}$ and $\{(x, t_x)\mid \hat{\tau}(x, t_x)\leq t_0, x\in U\}\setminus\{(x_0, t_0)\}$
is homotopy equivalent with $\{x\in U\mid f_{t_0}(x)\leq 0\}\setminus\{x_0\}$, and then (\ref{5.2}) follows by
definition and the homotopy invariance of homology. (\ref{5.2}) is also true in an equivariant setting, and then
we apply it to the functional $F_{t, K^\prime}|_S$, where $S$ is the unit sphere of $W^{1, 2}(\R/\Z; \R^{2n})$
and $F_{t, K^\prime}|_S$ is the restriction on $S$, we have
\be C_{S^1, k}(\hat{\tau}, S^1\cdot(x_0,t_0)) \cong C_{S^1, k}(F_{t_0, K^\prime}|_S, S^1\cdot x_0)=0,   \lb{5.3}\ee
where $k\in I$ and $t_0\in [a, a^\prime]$, $x_0$ corresponds to a $t_0 T$-periodic solution of (\ref{1.1}).
Then by Theorem 1.6.1 of \cite{Cha1} (which can be easily generalized to the equivariant sense), we obtain
(7.25) of \cite{Vit1}. The proof is complete. \hfill\hb

\setcounter{equation}{0}
\section{Proof of Theorem 1.1}

In this section, we give a proof for the Theorem 1.1 with $\wtd{H}_a(x)=a\vf_a(j(x))$, where $\vf_a$ satisfies
(i)-(ii) of Lemma 2.2.

Let $F_{a, K}$ be a functional defined by (\ref{2.6}) for some $a, K\in\R$ large enough and let $\vep>0$ be
small enough such that $[-\vep, 0)$ contains no critical values of $F_{a, K}$. We consider the exact sequence
of the triple $(X, X^{-\ep}, X^{-b})$ (for $b$ large enough)
\bea \rightarrow H_{S^1, *}(X^{-\ep}, X^{-b})
&\to& H_{S^1, *}(X,X^{-b})  \nn\\
&\to& H_{S^1, *}( X, X^{-\ep})\to H_{S^1, *-1}( X^{-\ep}, X^{-b})\to \cdots,  \lb{6.1}\eea
where $X=W^{1, 2}(\R/{T\Z}; \R^{2n})$. The normalized Morse series of $F_{a, K}$ in $ X^{-\ep}\setminus X^{-b}$
is defined, as usual, by
\be  M_a(t)=\sum_{q\ge 0,\;1\le j\le p} \dim C_{S^1,\;q}(F_{a, K}, \;S^1\cdot v_j)t^{q-d(K)},  \lb{6.2}\ee
where we denote by $\{S^1\cdot v_1, \ldots, S^1\cdot v_p\}$ the critical orbits of $F_{a, K}$ with critical
values less than $-\vep$. We denote by $t^{d(K)}H_a(t)$ the Poincar\'e series of $H_{S^1, *}(X^{-\ep}, X^{-b})$,
$H_a(t)$ is a Laurent series, and we have the equivariant Morse inequality
\be  M_a(t)-H_a(t)=(1+t)R_a(t),  \lb{6.3}\ee
where $R_a(t)$ is a Laurent series with nonnegative coefficients. On the other hand, the Poincar\'e series of
$H_{S^1, *}( X, X^{-b})$ is, by Corollary 5.11 of \cite{Vit1}, $t^{d(K)}(1/(1-t^2))$. The Poincar\'e series of
$H_{S^1, *}( X, X^{-\ep})$ is $t^{d(K)}Q_a(t)$, according to Theorem 5.1, if we set
$Q_a(t)=\sum_{k\in \Z}{q_kt^k}$, then
\be   q_k=0 \qquad\qquad \forall\;k\in \mathring {I},  \lb{6.4}\ee
where $I$ is defined in Theorem 5.1. Now by (\ref{6.1}) (cf. Proposition 1 in Appendix 2 of \cite{Vit1}), we
have
\be  H_a(t)-\frac{1}{1-t^2}+Q_a(t) = (1+t)S_a(t),   \lb{6.5}\ee
with $S_a(t)$ a Laurent series with nonnegative coefficients. Adding up (\ref{6.3}) and (\ref{6.5}) yields
\be  M_a(t)-\frac{1}{1-t^2}+Q_a(t) = (1+t)U_a(t),   \lb{6.6}\ee
where $U_a(t)=\sum_{i\in \Z}{u_it^i}$ also has nonnegative coefficients.

Now truncate (\ref{6.6}) at the degrees $2C$ and $2N$, where we set $C$ equal to $2n^2$, and $2N>2C$, and
write $M_a^{2N}(2C; t)$, $Q_a^{2N}(2C; t)\cdots$ for the truncated series. Then from (\ref{6.6}) we infer
\bea
&& M_a^{2N}(2C; t)-\sum_{h=C}^{N}{t^{2h}}+Q_a^{2N}(2C;t)   \nn\\
&& \qquad = (1+t)U_a^{2N-1}(2C; t)+t^{2N}u_{2N}+t^{2C}u_{2C-1}.  \lb{6.7}\eea
By (\ref{6.4}), and the fact that for $a$ large enough $\mathring {I}$ contains $[2C, 2N]$, indeed let $\alpha>0$ such that
any prime closed characteristic $(\tau, y)$ with $\hat{i}(y)\neq 0$ has $|\hat{i}(y)|>\alpha$. Then if
$k\geq aT/\min{\{\tau_i\}}$, we have $|i(y^k)| \sim k|\hat{i}(y)|\geq k\alpha\geq {a\alpha T}/\min{\{\tau_i\}}$,
which tends to infinity as $a\to\infty$, so $Q_a^{2N}(2C; t)=0$, and (\ref{6.7}) can be written
\be  M_a^{2N}(2C;t)-\sum_{h=C}^{N}{t^{2h}} = (1+t)U_a^{2N-1}(2C;t)+t^{2N}u_{2N}+t^{2C}u_{2C-1}.  \lb{6.8}\ee
Changing $C$ into $-C$, $N$ into $-N$, and counting terms with $-2N\leq i\leq -2C$, we obtain
\be  M_a^{-2C}(-2N; t) = (1+t)U_a^{-2C-1}(-2N;t)+t^{-2N}u_{-2N-1}+t^{-2C}u_{-2C}.  \lb{6.9}\ee
Denote by $\{x_1, \ldots, x_k\}$ the critical points of $F_{a, K}$ corresponding to $\{y_1,\ldots, y_k\}$. Note
that $v_1,\ldots,v_p$ in (\ref{6.2}) are iterations of $x_1,\ldots,x_k$. Since
$C_{S^1,\;q}(F_{a, K}, \;S^1\cdot x_j^m)$ can be non-zero only for $q=d(K)+i(y_j^m)+l$ with $0\le l\le 2n-2$,
by Propositions 2.11, 4.2 and Remark 4.4, the normalized Morse series (\ref{6.2}) becomes
\be  M_a(t) = \sum_{1\le j\le k,\; 0\le l\le 2n-2 \atop 1\le m_j<aT/\tau_j}k_l(y_j^{m_j})t^{i(y_j^{m_j})+l}
  \;\;= \sum_{1\le j\le k,\; 0\le l\le 2n-2 \atop 1 \le m_j\le K_j,\; sK_j+m_j<aT/\tau_j}
             k_l(y_j^{m_j})t^{i(y_j^{sK_j+m_j})+l}, \lb{6.10}\ee
where $K_j=K(y_j)$ and $s\in\N_0$. The last equality follows from Proposition 4.6.

Write $M(t)=\sum_{h\in \Z}w_ht^h$. Then we have
\be w_h\ = \sum_{1\le j\le k,\; 0\le l\le 2n-2 \atop 1 \le m\le K_j}
              k_l(y_j^m)\,^\#\{s\in\N_0\,|\,i(y_j^{sK_j+m})+l=h\}, \quad \forall\;2C\leq|h|\le 2N. \lb{6.11}\ee
Note that the right hand side of (\ref{6.10}) contains only those terms satisfying $sK_j+m_j<\frac{aT}{\tau_j}$.
Thus (\ref{6.11}) holds for $2C\leq |h|\leq 2N$ by (\ref{6.10}).

{\bf Claim 1.} {\it $w_h\le C_1$ for $2C\leq |h|\leq 2N$ with $C_1$ being independent of $a, K$}.

In fact, we have
\bea
^\#\{s\in\N_0 &|& i(y_j^{sK_j+m})+l=h \}  \nn\\
&=& \;^\#\{s\in\N_0 \;| \; i(y_j^{sK_j+m})+l=h,\;|i(y_j^{sK_j+m})-(sK_j+m)\hat{i}(y_j)|\le 2n\} \nn\\
&\le& \;^\#\{s\in\N_0 \;| \;|h-l-(sK_j+m)\hat{i}(y_j)|\le 2n\}  \nn\\
&=& \;^\#\left\{s\in\N_0 \; \left|\;\frac{}{}\right.\;h-l-2n-m\hat{i}(y_j)\le sK_j\hat{i}(y_j)
       \le h-l+2n-m\hat{i}(y_j)\right\}  \nn\\
&\le&\; \frac{4n}{K_j|\hat{i}(y_j)|}+2,  \lb{6.12}\eea
where the first equality follows from the fact
\be  |i(y_j^m)-m\hat{i}(y_j)|\le 2n,\quad \forall m\in\N,\; 1\le j\le k,  \lb{6.13}\ee
which follows from Theorems 10.1.2 of \cite{Lon4} and Theorem 2.1 of \cite{HuL1}, Note that
$i(y_j^{sK_j+m})+l=h\in[2C, 2N]$ holds only when $\hat{i}(y_j)>0$ and $i(y_j^{sK_j+m})+l=h\in[-2N, -2C]$ holds
only when $\hat{i}(y_j)< 0$. Hence Claim 1 holds.

Next we estimate $M_a^{2N}(2C; -1)$ and $M_a^{-2C}(-2N; -1)$. By (\ref{6.11}) we obtain
\bea
&& M_a^{2N}(2C; -1) = \sum_{h=2C}^{2N} w_h(-1)^h   \nn\\
&&\qquad = \sum_{1\le j\le k,\; 0\le l\le 2n-2 \atop 1 \le m\le K_j}(-1)^{i(y_j^m)+l}k_l(y_j^m)
              \,^\#\{s\in\N_0 \,|\, 2C\le i(y_j^{sK_j+m})+l\le 2N\}.  \lb{6.14}\eea
Here the second equality holds by (4.25). Similarly, we have
\bea
&& M_a^{-2C}(-2N; -1) = \sum_{h=-2N}^{-2C} w_h(-1)^h   \nn\\
&&\qquad = \sum_{1\le j\le k,\; 0\le l\le 2n-2 \atop 1 \le m\le K_j}(-1)^{i(y_j^m)+l}k_l(y_j^m)
              \,^\#\{s\in\N_0 \,|\, -2N\le i(y_j^{sK_j+m})+l\le -2C\}.  \qquad \lb{6.15}\eea

{\bf Claim 2.} {\it There is a real constant $C_2>0$ independent of $a,K$ such that
\bea
\left|M_a^{2N}(2C; -1)-\sum_{1\le j\le k,\;0\le l\le 2n-2 \atop 1 \le m\le K_j, \hat{i}(y_j)>0}
            (-1)^{i(y_j^m)+l}k_l(y_j^m)\frac{2N}{K_j\hat{i}(y_j)}\right|
               &\le& C_2,  \lb{6.16}\\
\left|M_a^{-2C}(-2N; -1)-\sum_{1\le j\le k,\; 0\le l\le 2n-2 \atop 1 \le m\le K_j, \hat{i}(y_j)<0}
            (-1)^{i(y_j^m)+l}k_l(y_j^m)\frac{2N}{K_j\hat{i}(y_j)}\right|
               &\le& C_2,  \lb{6.17}\eea
where the sum in the left hand side of $(6.16)$ equals to
$\;2N\sum_{\hat i(y_j)>0}\frac{\hat\chi(y_j)}{\hat i(y_j)}\;$, the sum in the left hand side of
(\ref{6.17}) equals to $\;2N\sum_{\hat i(y_j)<0}\frac{\hat\chi(y_j)}{\hat i(y_j)}\;$ by (\ref{4.32}).}

In fact, we have the estimates
\bea  ^\#\{s\in\N_0 &|& 2C\le i(y_j^{sK_j+m})+l\le 2N\}   \nn\\
&=& \;^\#\{s\in\N_0 \;| \; 2C\le i(y_j^{sK_j+m})+l\le 2N,\;
             |i(y_j^{sK_j+m})-(sK_j+m)\hat{i}(y_j)|\le 2n\}  \nn\\
&\le& \;^\#\{s\in\N_0 \;| \;0< (sK_j+m)\hat{i}(y_j)\le 2N-l+2n\}  \nn\\
&=& \;^\#\left\{s\in\N_0 \; \left |\;\frac{}{}\right.
     \;0\le s\le \frac{2N-l+2n-m\hat{i}(y_j)}{K_j\hat{i}(y_j)}\right\}  \nn\\
&\le& \; \frac{2N-l+2n}{K_j\hat{i}(y_j)}+1.  \nn\eea
On the other hand, we have
\bea
^\#\{s\in\N_0 &|& 2C\le i(y_j^{sK_j+m})+l\le 2N\}  \nn\\
&=& \;^\#\{s\in\N_0 \;| \; 2C\le i(y_j^{sK_j+m})+l\le 2N,\;
               |i(y_j^{sK_j+m})-(sK_j+m)\hat{i}(y_j)|\le 2n\}  \nn\\
&\ge& \;^\#\{s\in\N_0\;|\;i(y_j^{sK_j+m})\le(sK_j+m)\hat{i}(y_j)+2n\le 2N-l\} \nn\\
&\ge& \;^\#\left\{s\in\N_0 \; \left |\;\frac{}{}\right.
        \;0\le s\le \frac{2N-l-2n-m\hat{i}(y_j)}{K_j\hat{i}(y_j)}\right\}  \nn\\
&\ge& \;\frac{2N-l-2n}{K_j\hat{i}(y_j)}-2,  \nn\eea
where $m\le K_j$ is used and we note that $\hat{i}(y_j)> 0$ when $2C\le i(y_j^{sK_j+m})+l\le 2N$. Combining
these two estimates together with (\ref{6.14}), we obtain (\ref{6.16}). Similarly, we obtain (\ref{6.17}).

Note that all coefficients of $U_a(t)$ in (\ref{6.8}) and (\ref{6.9}) are nonnegative; hence, by Claim 1,
we have $u_{h}\leq w_h\leq C_1$ for $h=2N$ or $-2C$ and $u_{h}\leq w_{h+1}\leq C_1$ for $h=2C-1$ or $-2N-1$.
Now we choose $a$ to be sufficiently large, then we choose $N$ to be sufficiently large.

Note that by Claims 1 and 2, the constants $C_1$ and $C_2$ are independent of $a$ and $K$. Hence
dividing both sides of (\ref{6.8}), (\ref{6.9}) by $2N$ and letting $t=-1$, we obtain
\bea
M_a^{2N}(2C;-1)-(N-C+1) &=& u_{2N}+u_{2C-1}, \nn\\
M_a^{-2C}(-2N; -1) &=& u_{-2N-1}+u_{-2C}. \nn\eea
Dividing both sides of the above two identities by $2N$ and letting $N$ tend to infinity, we obtain
\bea
&&\lim_{N\to\infty}\frac{1}{2N}M_a^{2N}(2C; -1) = \frac{1}{2}, \nn\\
&&\lim_{N\to\infty}\frac{1}{2N}M_a^{-2C}(-2N; -1)= 0. \nn\eea
Hence (\ref{1.2}) and (\ref{1.3}) follow from (\ref{6.16}) and (\ref{6.17}).

Let us also mention that if there is no solution with $\hat{i}=0$, we do not need to cut our series at
$\pm 2C$; we can cut at $-2N$ and $2N$ only, thus obtaining
\be   M(t)-\frac{1}{1-t^2}=(1+t)U(t),   \lb{6.18}\ee
where $M(t)$ denotes $M_a(t)$ as $a$ tends to infinity. \hfill\hb

{\bf Acknowledgements.} The authors would like to sincerely thank the anonymous referee for his/her
careful reading of the manuscript and valuable comments. And we also would like to sincerely thank
Professor Erasmo Caponio for his valuable comments to us on the first manuscript arXiv:1308.3543v1
of this paper in 2013.

\bibliographystyle{abbrv}

\end{document}